\documentclass[12pt,twoside]{article}
\usepackage{amsmath,amsbsy,amsfonts}
\newtheorem{lemma}{Lemma}[section]
\newtheorem{proposition}{Proposition}[section]
\newtheorem{theorem}{Theorem}[section]
\newtheorem{corollary}{Corollary}[section]

\let\Section=\section
\def\section{\setcounter{equation}{0}\Section}

\begin{document}

\title{\Large\sf{Multiple solutions for a quasilinear Schr\"{o}dinger \\ equation on $\mathbb{R}^{N}$}}

\date{}

\author{\sf{Claudianor O. Alves} \thanks{Research partially supported by INCT-MAT, PROCAD and CNPq/Brazil 303080/2009-4}\\
{\small \textit{Unidade Acad\^emica de Matem\'atica e Estat\'{\i}stica}}\\
{\small \textit{Universidade Federal de Campina Grande}}\\
{\small \textit{58429-900, Campina Grande - PB - Brazil}}\\
{\small \textit{e-mail address: coalves@dme.ufcg.edu.br}} \\
\\
\vspace{1mm} \sf{Giovany M. Figueiredo} \thanks{Research partially
supported by supported by CNPq/Brazil
300705/2008-5.}\\
{\small \textit{Faculdade de Matem\'atica}}\\
{\small \textit{Universidade Federal do Par\'a}}\\
{\small \textit{66075-110, Bel\'em - PA - Brazil}}\\
{\small \textit{e-mail address: giovany@ufpa.br}}} \maketitle

%%% ----------------------------------------------------------------------
\begin{abstract}
\noindent The multiplicity of positive weak solutions for a quasilinear Schr\"{o}dinger equations $-L_p u
+(\lambda A(x)+1)|u|^{p-2}u= h(u)$ in $\mathbb{R}^N$ is
established, where $L_p u\doteq \epsilon^{p}\Delta_p u
+\epsilon^{p}\Delta_p (u^2)u$, $A$ is a nonnegative continuous
function and nonlinear term $h$ has a subcritical growth. We
achieved our results by using minimax methods and
Lusternik-Schnirelman theory of critical points.

\end{abstract}
\maketitle

\bigskip



\vskip.2cm

%------------------------------------------------------------------------------
\section{Introduction}
%------------------------------------------------------------------------------

In this paper we establish existence and
multiplicity of positive weak solutions for the following class of
quasilinear Schr\"{o}dinger equations:
\begin{equation}
\label{principal}\tag{$P_{\epsilon,\lambda}$} -L_pu +
(\lambda A(x)+1)|u|^{p-2}u = h(u),\ \ u\in W^{1,p}(\mathbb{R}^N),
\end{equation}
where
\[
L_p u\doteq \epsilon^{p}\Delta _pu +\epsilon^{p}\Delta _p(u^2)u,
\]
$\Delta_p u = div(|\nabla u |^{p-2} \nabla u)$ is the $p$-Laplacian, $\epsilon , \lambda $ are positive parameters, $2\leq p< N$ and function $A: \mathbb{R}^{N} \to \mathbb{R}$ satisfies the following conditions:
\begin{itemize}
\item[$(A_{1})$] \,\, $A \in C^{1}(\mathbb{R}^{N}, \mathbb{R}), A(x)
\geq 0$ for all $x \in \mathbb{R}^{N}$ and $\Omega={\rm{int}}
A^{-1}(0)$ is a nonempty bounded open set with smooth boundary
$\partial \Omega$ and $0 \in \Omega$. Moreover,
$A^{-1}(0)=\overline{\Omega} \cup D$ where $D$ is a set of measure
zero.

\item[$(A_{2})$] There exists $K_{0}>0$ such that
$$
\mu \Big(\Big\{ x \in \mathbb{R}^{N}: \, A(x) \leq K_{0}
\Big\}\Big)< \infty,
$$
\end{itemize}
where $\mu$ denotes the Lebesgue measure on $\mathbb{R}^{N}$.

On the nonlinearity $h$, we assume that it is of class $C^1$ and satisfies
the following conditions:
\begin{enumerate}
\item[$(H_1)$]
$h'(s)=o(|s|^{p-2})$ at the origin;
\item[$(H_2)$] $\displaystyle \lim_{|s|\rightarrow
\infty}h'(s)|s|^{-q+2}=0$ for some $q\in (2p,2p^{*})$ where
$p^*=Np/(N-p)$;
\item[$(H_3)$] There exists $\theta > 2p$ such that $0<\theta
H(s)\leq sh(s)$ for all $s>0$.
\item[$(H_4)$] The function $s\rightarrow h(s)/s^{2p-1}$ is
increasing for $s>0$.
\end{enumerate}

A typical example of a function satisfying the conditions
$(H_1)-(H_4)$ is given by $h(s)=s^{\mu}$ for $s\geq 0$, with
$2p-1<\mu<q-1$, and $h(s)=0$ for $s<0$.\\

For $p=2$, the solutions of (\ref{principal}) are related to existence of standing wave solutions for quasilinear
Schr\"{o}dinger equations of the form
\begin{equation}\label{SCH}
i\partial_t \psi =-\Delta
\psi+V(x)\psi-\widetilde{h}(|\psi|^2)\psi- \kappa\Delta
[\rho(|\psi|^2)]\rho'(|\psi|^2)\psi,
\end{equation}
where $\psi:\mathbb{R}\times  \mathbb{R}^N \rightarrow\mathbb{C}$,
$V$ is a given potential, $\kappa$ is a real constant and
$\rho,\widetilde{h}$ are real functions.  Quasilinear equations of
the form (\ref{SCH}) have been studied in relation with some
mathematical models in physics. For example, when $\rho(s)=s$, the
above equation is
\begin{equation}\label{SCH2}
i\partial_t \psi =-\Delta \psi+V(x)\psi- \kappa\Delta
[|\psi|^2]\psi-\widetilde{h}(|\psi|^2)\psi.
\end{equation}
It was shown that a system describing the self-trapped electron on
a lattice can be reduced in the continuum limit to (\ref{SCH2})
and numerics results on this equation are obtained in
\cite{Brizhik}. In \cite{Hartmann}, motivated by the nanotubes and
fullerene related structures, it was proposed and shown that a
discrete system describing the interaction of a 2-dimensional
hexagonal lattice with an excitation caused by an excess electron
can be reduced to (\ref{SCH2}) and numerics results have been done
on domains of disc type, cylinder type and sphere type. The
superfluid film equation in plasma physics  has also the structure
(\ref{SCH}) for $\rho(s) = s$, see \cite{Kurihura}.

The general equation (\ref{SCH}) with various form of quasilinear
terms $\rho(s)$ has been derived as models of several other
physical phenomena corresponding to various types of $\rho(s)$.
For example, in the case $\rho(s) = (1+s)^{1/2}$ , equation
(\ref{SCH}) models the self-channeling of a high-power ultra short
laser in matter, see \cite{Borovskii-Galkin} and \cite{Ritchie}.
Equation (\ref{SCH}) also appears in fluid mechanics
\cite{Kosevich-Ivanov}, in the theory of Heisenberg ferromagnets
and magnons \cite{Takeno-Homma}, in dissipative quantum mechanics
and in condensed matter theory \cite{Makhankov-Fedyanin}. The
Semilinear case corresponding to $\kappa=0$ in the whole
$\mathbb{R}^{N}$ has been studied extensively in recent years, see
for example \cite{Floer-Weinstein},
\cite{Jeanjean-Tanaka} and references therein.

Putting $\psi(t,x)=\exp(-iFt)u(x),\; F\in \mathbb{R}$, into the
equation (\ref{SCH2}), we obtain a corresponding equation
\begin{equation}\label{simp}
-\Delta u -\Delta (u^2)u +V(x)u= h(u)
\end{equation}
where we have renamed $V(x)-F$ to be $V(x)$,
$h(u)=\widetilde{h}(u^2)u$ and we assume, without loss of
generality, that $\kappa =1$.

The quasilinear equation (\ref{simp}) in the whole $\mathbb{R}^{N}$
has received special attention in the past several years, see for
example the works \cite{AlvesFigueiredoSevero}, \cite{AMS1},
\cite{AMS2}, \cite{MJAlves-Carriao-Miyagaki}, \cite{Jeanjean-Colin},
\cite{Jeanjean-Colin-Squassina}, \cite{do O-Severo}, \cite{OMS},
\cite{Liu-Wang I}, \cite{Liu-Wang II}, \cite{Liu-wang-wang},
\cite{Poppenberg-Schmitt-Wang}, \cite{Elves1}, \cite{Elves2} and
references therein. In these papers, we find important results on
the existence of nontrivial solutions of (\ref{simp}) and a good
insight into this quasilinear Schr\"{o}dinger equation. The main
strategies used are the following: the first of them consists in by
using a constrained minimization argument, which gives a solution of
(\ref{simp}) with an unknown Lagrange multiplier $\lambda$ in front
of the nonlinear term, see for example
\cite{Poppenberg-Schmitt-Wang}. The other one consists in by using a
special change of variables to get a new semilinear equation and an
appropriate Orlicz space framework, for more details see
\cite{Jeanjean-Colin}, \cite{do O-Severo} and \cite{Liu-Wang II}. In
\cite{AlvesFigueiredoSevero}, existence, multiplicity and
concentration of solutions have been study, by using the
Lusternik-Schnirelman category, for the following class of problems
$$
-\epsilon^{p}\Delta _pu -\epsilon^{p}\Delta _p(u^2)u+
V(x)|u|^{p-2}u = h(u),\ \ u\in W^{1,p}(\mathbb{R}^N),
$$
when $\epsilon$ is sufficiently small and $V$ satisfying the condition
$$
\displaystyle\liminf_{|x|\rightarrow \infty}V(x) >
\displaystyle\inf_{x \in \mathbb{R}^{N}} V(x)=V_0>0, \eqno{(R)}
$$
which has been introduced by Rabinowitz \cite{Rabinowitz}. In \cite{AlvesFigueiredoSevero1}, multiplicity of solutions also have been proved for problem of the type
$$
-\epsilon^{p}\Delta _pu -\epsilon^{p}\Delta _p(u^2)u = h(u),\ \ u\in W_{0}^{1,p}(\lambda \Omega),
$$
where $\Omega$ is a bounded domain and $\lambda$ is large enough.

Related to the $p$-Laplacian operator, we would like to cite a paper due to Alves and Soares \cite{AlvesSoares},   where existence, multiple and concentration of solutions, by using the Lusternik-Schnirelman category, have been established for the following class of $p$-Laplacian equations
$$
-\epsilon^{p}\Delta_{p}u + (\lambda A(x)+1)|u|^{p-2}u  = h(u),\ \ u\in W^{1,p}(\mathbb{R}^{N}), \eqno{(P)}
$$
by assuming that $A$ verifies conditions $(A_1)-(A_2)$, $h$ is a continuous functions with subcritical growth, $\epsilon$ is sufficiently small and $\lambda$ is large enough. In that paper, it is proved  that there exists $\epsilon^{*}>0$ such that for any $\epsilon \in (0,\epsilon^{*})$ there exists $\lambda^{*}(\epsilon)>0$ such that $(P)$ has at least $cat(\Omega)$ solutions for
any $\lambda \geq \lambda^{*}(\epsilon)$. The results showed in \cite{AlvesSoares} are associated with the main results proved by  Barstch and Wang \cite{BW2,BW1}, where condition $(A_2)$ was used by the first time. Here, we would like to emphasize that the assumptions $(A_{1})-(A_{2})$ do
not imply that potential $V(x)=\lambda A(x)+1$ verifies condition $(R)$.

The present paper was motivated by works \cite{AlvesFigueiredoSevero}, \cite{AlvesFigueiredoSevero1} and \cite{AlvesSoares}. Here, we intend to show that the same type of results found in \cite{AlvesSoares} for $p$-Laplacian also hold for operator $L_p$. However, due to the presence of the term $\Delta_{p}(u^{2})u$ in  $L_{p}u$, several estimates used \cite{AlvesSoares} can not be repeated  for the functional energy associated to (\ref{principal}), given by
\[
J_{\epsilon, \lambda}(u) = \frac{1}{p}\int_{\mathbb{R}^N}
\epsilon^{p}(1+2^{p-1}|u|^p)|\nabla u |^p  +
\frac{1}{p}\int_{\mathbb{R}^N} (\lambda A(x) + 1)|u|^p  -
 \int_{\mathbb{R}^N} H(u),
\]
where $H(s)=\int_0^s h(t)dt$. As observed in \cite{severo} and \cite{severo2} , there are some technical difficulties to apply directly variational methods to $J_{\epsilon, \lambda}$.  The main difficult is related to the fact that $J_{\epsilon, \lambda}$ is not well defined in $W^{1,p}(\mathbb{R}^N)$. By a direct computation, if $u\in C_0^1(\mathbb{R}^N\backslash \{0\})$ is defined by
\[
u(x)=|x|^{(p-N)/2p}\ \ \mbox{for}\ \  x\in B_1\backslash \{0\},
\]
then $u\in W^{1,p}(\mathbb{R}^N)$ while the function $|u|^p|\nabla
u |^p$ does not belong to $L^1(\mathbb{R}^N)$. To overcome this difficulty, we use a change variable developed in
\cite{severo} and \cite{severo2}, which generalize one found in Liu, Wang and Wang \cite{Liu-Wang II} and Colin-Jeanjean
\cite{Jeanjean-Colin} for the case $p=2$.

Before to state our main result, we recall that if $Y$ is a closed set of a topological space $X$, we denote the Lusternik-Schnirelman category of $Y$ in $X$ by $cat_X (Y)$, which is the least number of closed and contractible sets in $X$ that cover $Y$. Hereafter, $cat \, X$ denotes $cat_X (X)$.

\vspace{0.5 cm}

 The main result that we prove is the following:

\begin{theorem} \label{T1}
Suppose that $(A_{1})-(A_{2})$ and $(H_{1})-(H_{4})$ hold.  Then
there exists $\epsilon^{*}>0$ such that for any $\epsilon \in (0,
\epsilon^{*})$ there exists $\lambda^{*}(\epsilon)>0$ such that
$(P_{\epsilon,\lambda})$ has at least $cat(\Omega)$ solutions
for any $\lambda \geq \lambda^{*}(\epsilon)$.
\end{theorem}

To finish this introduction, we would like to emphasize that Theorem
\ref{T1} can be seen as a complement of the studies made in
\cite{AlvesFigueiredoSevero} and \cite{AlvesSoares}. In \cite{AlvesFigueiredoSevero}, the authors considered a class
of problems involving the $L_p$ operator, however the
potential $V$ verifies the condition $(R)$, while that in \cite{AlvesSoares}, only the $p$-Laplacian was considered.

The plan of this paper is as follows. In Section 2, we review some proprieties of the change variable that we will apply.  Section 3 establishes a compactness result for the energy functional for all sufficiently large  $\lambda$  and arbitrary $\epsilon$. In Section 4 is made a
study of the behavior of the minimax levels with respect to parameter $\lambda$ and $\epsilon$. Section 5  offers the proof of our main result.

\section{Variational framework and preliminary results}

Since we intend to find positive solutions, let us assume that
$$
h(s)=0 \ \ \mbox{for all} \ \ s<0.
$$

Moreover, hereafter, we will work with the following problem
equivalent to $(P_{\epsilon,\lambda})$, which is obtained under
change of variable $\epsilon z=x$
$$
\ \  \left\{
             \begin{array}{l}
              - \Delta_{p}u - \Delta_{p}(u^{2})u + (\lambda
A(\epsilon x)+1)|u|^{p-2}u=h(u)
               \ \ in \ \ \mathbb{R}^N
        \\
        u \in W^{1,p}(\mathbb{R}^N) \,\,\, \mbox{with} \,\,\, 2\leq p< N,
        \\
        u(x) > 0, \forall \  x \in \mathbb{R}^N .
             \end{array}
           \right.
           \eqno{(P^{*}_{\epsilon,\lambda})}
$$
In what follows, we use the change variable developed in
\cite{severo} and \cite{severo2} which generalizes one found in Liu, Wang and Wang \cite{Liu-Wang II} and Colin-Jeanjean
\cite{Jeanjean-Colin} for the case $p=2$. More precisely, let us introduce the change of variables $v= f^{-1}(u)$,
where $f$ is defined by
\begin{equation}\label{mudanda de variavel}
\begin{array}{cllc}
f'(t)&=&\dfrac{1}{(1+ 2^{p-1}|f(t)|^p)^{1/p}} & \;\;\mbox{ on }\;\;  [0,+\infty),\\
 f(t)&=&-f(-t) & \;\;\mbox{ on }\;\; (-\infty,0].
\end{array}
\end{equation}
Therefore, using the above change of variables, we consider a new functional
$I_{\epsilon,\lambda}$, given by
\begin{equation}\label{funcional}
I_{\epsilon,\lambda}(v)\doteq J_{\epsilon,\lambda}(f(v))=
\frac{1}{p}\int_{\mathbb{R}^N} |\nabla v |^p\   +
\frac{1}{p}\int_{\mathbb{R}^N} (\lambda A(\epsilon x)+1)|f(v)|^p\
  -
 \int_{\mathbb{R}^N} H(f(v))\
\end{equation}
which is well defined on the Banach space $X$ defined
by
$$
X_{\epsilon,\lambda} = \left\{ u \in W^{1 , p}({\mathbb{R}^N}) :
\displaystyle\int_{\mathbb{R}^N}A(\epsilon x)|f(v)|^{p} < \infty
\right\}
$$
endowed with the norm
$$
\|u\|_{{\epsilon,\lambda}} = |\nabla u|_{p} +
\inf_{\xi>0}\frac{1}{\xi}\biggl[1+\displaystyle\int_{\mathbb{R}^N}(\lambda
A(\epsilon x)+1)|f(\xi v)|^{p}\biggl].
$$

A direct computation shows that $I_{\epsilon,\lambda}:X_{\epsilon,\lambda} \rightarrow \mathbb{R}$
is of class $C^1$ under the conditions $(A_1)-(A_2)$ and
$(H_1)-(H_2)$ with
$$
\begin{array}{l}
I_{\epsilon,\lambda}'(v)w= \displaystyle \int_{\mathbb{R}^N}|\nabla v|^{p-2}\nabla
v\nabla w\ +(\lambda A(\epsilon
x)+1)|f(v)|^{p-2}f(v)f'(v)w- \\
\mbox{} \\
\hspace{2 cm} - \displaystyle \int_{\mathbb{R}^N}h(f(v))f'(v)w
\end{array}
$$
for $v,w\in W^{1,p}(\mathbb{R}^N)$. Thus, the critical points of
$I_{\epsilon,\lambda}$ correspond exactly to the weak solutions of the problem
$$
\ \  \left\{
             \begin{array}{l}
              -\Delta_{p}v + (\lambda
A(\epsilon x)+1)|f(v)|^{p-2}f(v)f'(v)=h(f(v))f'(v)
               \ \ \mbox{in} \ \ \mathbb{R}^N
        \\
        v \in W^{1,p}(\mathbb{R}^N) \,\,\, \mbox{with} \,\,\, 2\leq p< N,
        \\
        v(x) > 0, \forall \  x \in \mathbb{R}^N .
             \end{array}
           \right.
           \eqno{(S_{\epsilon,\lambda})}
$$

\vspace{0.5 cm}

The below result establishes a relation between the solutions of $(S_{\epsilon,\lambda})$ with one of
$(P^{*}_{\epsilon,\lambda})$:

\begin{proposition} \label{prop1} If $v \in W^{1,p}(\mathbb{R}^N)\cap L^{\infty}_{loc}(\mathbb{R}^N)$ is a critical point
of the functional $I_{\epsilon,\lambda}$, then $u=f(v)$ is a weak solution of $(P^{*}_{\epsilon,\lambda})$.
\end{proposition}
\noindent {\bf Proof.} \, See \cite{severo2}. \hfill\rule{2mm}{2mm}

\medskip

From the above proposition, it is clear that to obtain a weak solution of
$(P_{\epsilon,\lambda})$, it is sufficient to obtain a critical
point of the functional $I_{\epsilon,\lambda}$ in
$W^{1,p}(\mathbb{R}^N)\cap L^{\infty}_{loc}(\mathbb{R}^N)$.

\vspace{0.5 cm}

Next, let us collect some properties of the change of variables
$f:\mathbb{R}\rightarrow \mathbb{R}$ defined in (\ref{mudanda de
variavel}), which will be usual in the sequel of the paper.

\begin{lemma}
\label{Lema f} The function $f(t)$ and its derivative enjoy the
following properties:
\begin{itemize}
  \item[(1)] $f$ is uniquely defined, $C^2$ and invertible;
  \item[(2)] $|f'(t)|\leq 1$ for all $t\in \mathbb{R}$;
  \item[(3)] $|f(t)|\leq |t|$ for all $t\in \mathbb{R}$;
  \item[(4)] $f(t)/t\rightarrow 1$ as $t\rightarrow
0$;
  \item[(5)] $|f(t)|\leq 2^{1/2p}|t|^{1/2}$ for all $t\in
  \mathbb{R}$;
  \item[(6)] $f(t)/2\leq tf'(t)\leq f(t)$ for all $t\geq 0$;
  \item[(7)] $f(t)/{\sqrt t}\rightarrow a>0$ as $t\rightarrow +\infty$.
  \item[(8)] there exists a positive constant
$C$ such that
\[
|f(t)| \geq
\begin{cases}
C|t|,\quad & |t| \leq 1 \\
C|t|^{1/2},\quad & |t|  \geq 1.
\end{cases}
\]
\item[$(9)$] $|f(t)f'(t)|\leq 1/2^{(p-1)/p}$ for all $t\in \mathbb{R}$.
\end{itemize}
\end{lemma}

\noindent {\bf Proof.} \, See \cite{severo2}. \hfill\rule{2mm}{2mm}

\vspace{0.5 cm}

The next lemma can be found in \cite{AlvesFigueiredoSevero}, however for convenience of the reader we will write its proof.

\begin{lemma}\label{resultadoClaudianor}
Let $(v_{n})$ be a sequence in
$W^{1,p}(\mathbb{R}^{N})$ verifying $\displaystyle\int_{\mathbb{R}^{N}}
|f(v_{n})|^{p}\rightarrow 0$ as $n \to \infty$. Then,
$$
\displaystyle\inf_{\xi
>0}\frac{1}{\xi}\left\{1+\displaystyle\int_{\mathbb{R}^{N}}|f(\xi
v_{n})|^{p}\right\}\rightarrow 0 \,\,\, \mbox{as} \,\,\, n \to
\infty.
$$
\end{lemma}
\noindent {\bf Proof.} Hereafter, once that $f$ is odd,  we can assume without loss of generality that $v_n \geq 0$ for all $n \in \mathbb{N}$. Since $f(t)/t$ is nonincreasing for $t>0$, for each $\xi>1$,
$$
\frac{1}{\xi}+ \frac{1}{\xi}\displaystyle\int_{\mathbb{R}^{N}}|f(\xi v_{n})|^{p}\leq \frac{1}{\xi}+\xi^{p-1}\int_{\mathbb{R}^{N}}|f( v_{n})|^{p}.
$$
Hence, for each $\delta >0$, fixing $\xi_{*}$ sufficiently large such that $\frac{1}{\xi_{*}}< \frac{\delta}{2}$, we get
$$
\displaystyle\inf_{\xi >0}\frac{1}{\xi}\left\{1 +\int_{\mathbb{R}^{N}}|f(\xi
v_{n})|^{p}\right\}\leq \frac{\delta}{2}+
\xi^{p-1}_{*}\int_{\mathbb{R}^{N}}|f(v_{n})|^{p}.
$$
Thus,
$$
\limsup_{n \to \infty}\left(\displaystyle\inf_{\xi
>0}\frac{1}{\xi}\biggl\{1+\int_{\mathbb{R}^{N}}|f(\xi
v_{n})|^{p}\biggl\}\right)\leq \frac{\delta}{2} \,\,\,\,\,\,
\textrm{for all } \delta >0
$$
which proves the proposition. \hfill\rule{2mm}{2mm}

\vspace{0.5 cm}

Repeating the same type of arguments explored in the proof of the last lemma, we have the following result

\begin{corollary}\label{resultadoClaudianor1}
Let $(v_n)$ be a sequence in $ X_{\epsilon,\lambda} $ with
$$\displaystyle\int_{\mathbb{R}^N} (\lambda A(\epsilon
x)+1)|f(v_{n})|^{p}\rightarrow 0 \,\, \mbox{as} \,\, n \to +\infty.
$$
Then
$$\displaystyle\inf_{\xi
>0}\frac{1}{\xi}\biggl\{1+\displaystyle\int_{\mathbb{R}^N}
(\lambda A(\epsilon x)+1)|f(\xi v_{n})|^{p}\biggl\}\rightarrow 0 \,\, n \to +\infty.
$$
\end{corollary}

\vspace{0.5 cm} The next lemma is related to a claim made in \cite{AlvesFigueiredoSevero}, which wasn't proved in that paper. Here, we decide to show its proof.

\begin{lemma} \label{Normaequivalente}
The function
$$
\|v\|=|\nabla v|_{p} +
\inf_{\xi>0}\frac{1}{\xi}\biggl[1+\displaystyle\int_{\mathbb{R}^N}|f(\xi v)|^{p}\biggl].
$$
is a norm in $W^{1,p}(\mathbb{R}^{N})$. Moreover, $\|\,\,\,\,\|$ is equivalent to the usual norm in  $W^{1,p}(\mathbb{R}^{N})$.
\end{lemma}

\noindent {\bf Proof.}\, We will omit the proof that $\|\,\,\,\|$ is a norm, because we can repeat with few modifications, the same arguments used by Severo \cite{severo}. From the hypotheses on $f$,
$$
0 \leq |f(t)| \leq |t| \,\,\,\, \forall t \in \mathbb{R}^{N},
$$
this way
$$
\int_{\mathbb{R}^{N}}|f(\xi v)|^{p} \leq \xi^{p} \int_{\mathbb{R}^{N}}|v|^{p} \,\,\, \forall \xi \geq 0,
$$
from where it follows that
$$
\inf_{\xi >0} \frac{1}{\xi}\left\{1+\int_{\mathbb{R}^{N}}|f(\xi v)|^{p} \right\} \leq \inf_{\xi >0} \left\{\frac{1}{\xi}+L \xi^{p-1} \right\}
$$
where
$$
L=\int_{\mathbb{R}^{N}}|v|^{p}.
$$
Now, let us consider the function
$$
g(\xi)=\frac{1}{\xi}+L\xi^{p-1}\,\,\,\,\, \mbox{for} \,\,\,\ \xi >0.
$$
A direct computation implies that $g$ has a global minimum at some $\xi_0 >0$, which satisfies
$$
g'(\xi_0)=0 \Leftrightarrow - {\xi_0}^{2}+(p-1)L{\xi_0}^{p-2}=0.
$$
Then,
$$
\xi_0 = \left( \frac{1}{(p-1)L} \right)^{\frac{1}{p}}
$$
and so,
$$
g(\xi_0)=(L(p-1))^{\frac{1}{p}}+L\left(\frac{1}{(p-1)L}\right)^{\frac{p-1}{p}}=CL^{\frac{1}{p}}
$$
for some $C>0$. Using these informations, it follows that
$$
\|v\| \leq |\nabla v|_{p}+C|v|_{p} \,\,\, \forall v \in W^{1,p}(\mathbb{R}^{N}).
$$
Hence, there is $c_1 >0$ such that
$$
\|v\| \leq  c_1 \|v\|_{1,p} \,\,\, \forall v \in W^{1,p}(\mathbb{R}^{N}).
$$
where
$$
\|v\|_{1,p}= \left[\int_{\mathbb{R}^{N}}|\nabla v|^{p}+\int_{\mathbb{R}^{N}}|v|^{p}\right]^{\frac{1}{p}} \,\,\, \forall v \in W^{1,p}(\mathbb{R}^{N}).
$$
Since $(W^{1,p}(\mathbb{R}^{N}), \|\,\,\,\, \|)$ and $(W^{1,p}(\mathbb{R}^{N}), \|\,\,\,\, \|_{1,p})$ are Banach spaces, the last inequality together with Closed Graphic Theorem yields  $\|\,\,\,\|$ and $\|\,\,\,\|_{1,p}$ are equivalent norms. \hfill\rule{2mm}{2mm}

\begin{lemma} \label{CZ}
Let $(v_{n})$ be a  sequence in  $W^{1,p}(\mathbb{R}^{N})$ and set
$$
\mathcal{Q}(v):=\int_{\mathbb{R}^{N}}|\nabla
v|^p +\int_{\mathbb{R}^{N}}
|f(v)|^{p},
$$
and
$$
\|v\|=|\nabla v|_{p} +
\inf_{\xi>0}\frac{1}{\xi}\biggl[1+\displaystyle\int_{\mathbb{R}^{N}}|f(\xi v)|^{p}\biggl].
$$
Then, $\mathcal{Q}(v_{n})\rightarrow 0$ if, and only
if, $\|v_{n}\|\rightarrow 0$. Moreover, $(v_n)$ is bounded in $(W^{1,p}(\mathbb{R}^{N}),\|\,\,\,\|)$ if, and only if, $(\mathcal{Q}(v_n))$ is bounded in $\mathbb{R}$.
\end{lemma}

\noindent {\bf Proof.} The first part of the lemma is an immediate consequence of Lemma \ref{resultadoClaudianor}, this way we will prove only the second part of the lemma.

A straightforward computation gives
$$
\|v \| \leq \mathcal{Q}(v)\,\,\, \forall v \in W^{1,p}(\mathbb{R}^{N}),
$$
from where it follows that if $(\mathcal{Q}(v_n))$ is bounded, then $(v_n)$  is also bounded. On the other hand, by Lemma \ref{Normaequivalente}, $(v_n)$ is a bounded sequence in  $(W^{1,p}(\mathbb{R}^{N}), \|\,\,\, \|)$ if, and only if, $(v_n)$ is bounded $(W^{1,p}(\mathbb{R}^{N}), \|\,\,\, \|_{1,p})$, where $\|\,\,\, \|_{1,p}$ is the usual norm in  $W^{1,p}(\mathbb{R}^{N})$. Hence, there is $M>0$ such that
$$
\int_{\mathbb{R}^{N}}|\nabla v_n|^{p} \leq M \,\,\, \mbox{and} \,\,\, \int_{\mathbb{R}^{N}}|v_n|^{p} \,\,\, \forall n \in \mathbb{N}.
$$
Recalling that
$$
|f(t)| \leq |t| \,\,\,\, \forall t \geq 0,
$$
we have the estimate
$$
\int_{\mathbb{R}^{N}}|f(v_n)|^{p} \leq \int_{\mathbb{R}^{N}}|v_n|^{p} \leq M \,\,\, \forall n \in \mathbb{N},
$$
which shows that $(\mathcal{Q}(v_n))$ is bounded.  \hfill\rule{2mm}{2mm}

\begin{lemma} \label{Convexa0} The function $|f|^p$ is a convex function, and so,
$$
(|f(t)|^{p-2}f(t)f'(t) - |f(s)|^{p-2}f(s)f'(s))(t-s) \geq 0 \,\,\,\,\, \forall t,s \in \mathbb{R}.
$$
\end{lemma}

\noindent {\bf Proof.} A direct computation shows that second derivative of the function
$$
Q(t)=|f(t)|^{p} \,\,\, \mbox{for} \,\, t \in \mathbb{R}
$$
satisfies the equality
$$
Q''(t)=\frac{p|f(t)|^{p-2}|f'(t)|^{2}\left((p-1)+(p-2)2^{p-1}|f(t)|^{p}\right)}{1+2^{p-1}|f(t)|^{p}} > 0 \,\,\, \forall t \in \mathbb{R} \setminus \{0\},
$$
implying that $Q$ is a convex function. From this,
$$
(Q'(t)-Q'(s))(t-s) \geq 0 \,\,\, \forall t,s \in \mathbb{R}
$$
that is,
$$
(|f(t)|^{p-2}f(t)f'(t) - |f(s)|^{p-2}f(s)f'(s))(t-s) \geq 0 \,\,\,\,\, \forall t,s \in \mathbb{R},
$$
finishing the proof.  \hfill\rule{2mm}{2mm}

\begin{lemma}\label{Convexa} Let $(v_n) \subset W^{1,p}(\mathbb{R}^{N})$ be a sequence of nonnegative functions such that $v_n \rightharpoonup v $ in $W^{1,p}(\mathbb{R}^{N})$, $v_n(x) \to v(x)$ a.e in $\mathbb{R}^{N}$ and
$$
\int_{\mathbb{R}^{N}}(|f(v_n)|^{p-2}f(v_n)f'(v_n)-|f(v)|^{p-2}f(v)f'(v))(v_n-v) \to 0 \,\,\, \mbox{as}  \,\,\, n \to +\infty.
$$
Then,
$$
\int_{\mathbb{R}^{N}}|f(v_n-v)|^{p} \to 0 \,\,\, \mbox{as}  \,\,\, n \to +\infty.
$$
\end{lemma}

\noindent {\bf Proof.}
By hypothesis,
$$
\int_{\mathbb{R}^{N}}(|f(v_n)|^{p-2}f(v_n)f'(v_n)-|f(v)|^{p-2}f(v)f'(v))(v_n-v)=o_n(1)
$$
or equivalently,
$$
\begin{array}{l}
\displaystyle \int_{\mathbb{R}^{N}}|f(v_n)|^{p-2}f(v_n)f'(v_n)v_n = \int_{\mathbb{R}^{N}}|f(v_n)|^{p-2}f(v_n)f'(v_n)v +  \\
\mbox{}\\
 \hspace{5 cm} \displaystyle \int_{\mathbb{R}^{N}}|f(v)|^{p-2}f(v)f'(v)(v_n-v) + o_n(1).
\end{array}
$$
Once that $v_n \rightharpoonup v$ in $W^{1,p}(\mathbb{R}^{N})$,
$$
\int_{\mathbb{R}^{N}}|f(v)|^{p-2}f(v)f'(v)(v_n-v)=o_n(1)
$$
and so,
$$
\int_{\mathbb{R}^{N}}|f(v_n)|^{p-2}f(v_n)f'(v_n)v_n = \int_{\mathbb{R}^{N}}|f(v_n)|^{p-2}f(v_n)f'(v_n)v + o_n(1).
$$
Recalling that
$$
|f(t)| \leq |t|  \,\,\,\, \mbox{and} \,\,\, |f'(t)| \leq 1 \,\,\, \forall t \in \mathbb{R},
$$
it follows that $(|f(v_n)|^{p-2}f(v_n)f'(v_n))$ is bounded sequence in $L^{\frac{p}{p-1}}(\mathbb{R}^{N})$. Hence,
$$
\int_{\mathbb{R}^{N}}|f(v_n)|^{p-2}f(v_n)f'(v_n)v \to \int_{\mathbb{R}^{N}}|f(v)|^{p-2}f(v)f'(v)v
$$
which gives
$$
\int_{\mathbb{R}^{N}}|f(v_n)|^{p-2}f(v_n)f'(v_n)v_n \to \int_{\mathbb{R}^{N}}|f(v)|^{p-2}f(v)f'(v)v.
$$
From Lemma \ref{Lema f},
$$
|f(t)|^{p} \leq 2|f(t)|^{p-2}f(t)f'(t)t \,\,\,\, \forall t \geq 0
$$
then,
$$
|f(v_n)|^{p} \leq |f(v_n)|^{p-2}f(v_n)f'(v_n)v_n  \,\,\,\, \forall n \in \mathbb{N}.
$$
Using the above informations together with Lebesgue's Theorem, we deduce
$$
\int_{\mathbb{R}}|f(v_n)|^{p} \to \int_{\mathbb{R}}|f(v)|^{p}.
$$
On the other hand, since $|f'(t)| \leq 1$ for all $t \in \mathbb{R}$, we have the inequality
$$
|f(v_n-v)|=f(|v_n-v|) \leq f(v_n+v) \leq f(v_n) +v \,\,\, \forall n \in \mathbb{N}
$$
which gives
$$
|f(v_n-v)|^{p} \leq 2^{p}(|f(v_n)|^{p} +|v|^{p} ) \,\,\, \forall n \in \mathbb{N}.
$$
Combining the last inequality  with Lebesgue's Theorem, we get
$$
\int_{\mathbb{R}^{N}}|f(v_n-v)|^{p} \to 0,
$$
concluding the proof of the lemma. \hfill\rule{2mm}{2mm}

\begin{corollary} \label{monotonia}
Let $(v_n) \subset W^{1,p}(\mathbb{R}^{N})$ be a sequence of nonnegative functions such that $v_n \rightharpoonup v $ in $W^{1,p}(\mathbb{R}^{N})$, $v_n(x) \to v(x)$ a.e in $\mathbb{R}^{N}$ and the below limits hold
\begin{equation} \label{Cor1}
\int_{\mathbb{R}^{N}}(|f(v_n)|^{p-2}f(v_n)f'(v_n)-|f(v)|^{p-2}f(v)f'(v))(v_n-v) \to 0 \,\,\, \mbox{as}  \,\,\, n \to +\infty
\end{equation}
and
\begin{equation} \label{Cor2}
\int_{\mathbb{R}^{N}}\left\langle |\nabla v_n|^{p-2}\nabla v_n - |\nabla v|^{p-2}\nabla v, \nabla v_n - \nabla v \right\rangle \to 0 \,\,\, \mbox{as}  \,\,\, n \to +\infty.
\end{equation}
Then, $v_n \to v$ in $W^{1,p}(\mathbb{R}^{N})$.
\end{corollary}

\noindent {\bf Proof.} By Lemma \ref{Convexa}, the limit (\ref{Cor1}) leads to
$$
\int_{\mathbb{R}^{N}}|f(v_n-v)|^{p} \to 0 \,\,\, \mbox{as}  \,\,\, n \to +\infty.
$$
On the other hand, the limit (\ref{Cor2}) implies that
$$
\int_{\mathbb{R}^{N}}|\nabla v_n - \nabla v|^{p} \to 0 \,\,\, \mbox{as}  \,\,\, n \to +\infty.
$$
The above limits give
$$
\mathcal{Q}(v_n - v) \to 0 \,\,\, \mbox{as}  \,\,\, n \to +\infty,
$$
and so, by Lemma \ref{CZ}
$$
\|v_n - v\| \to 0 \,\,\, \mbox{as}  \,\,\, n \to +\infty
$$
or equivalently,
$$
v_n \to v \,\,\, \mbox{in} \,\,\, W^{1,p}(\mathbb{R}^{N}),
$$
proving the lemma.  \hfill\rule{2mm}{2mm}

\vspace{0.5 cm}

%------------------------------------------------------------------------------
\section{The Palais-Smale condition}
%------------------------------------------------------------------------------

In this Section, the main goal is to show that $I_{\epsilon, \lambda}$
satisfies the Palais-Smale condition. To this end, we have to prove some
technical lemmas.

\begin{lemma}
\label{L1} Suppose that $h$ satisfies $(H_{1})-(H_{3})$. Let
$(v_{n}) \subset X_{\epsilon,\lambda}$ be a $(PS)_c$ sequence for
$I_{\epsilon,\lambda}$. Then there exists a constant $K>0$,
independent of $\epsilon$ and $\lambda$, such that
$$
\limsup_{n \to \infty}\|v_{n}\|_{\epsilon,\lambda} \leq K
$$
for all $\epsilon, \lambda>0$
\end{lemma}
\noindent {\bf Proof.} \, Using $(H_3)$ and  Lemma \ref{Lema f}$(6)$,
$$
c+o_{n}(1)\|v_{n}\|_{\epsilon,\lambda}\geq
(\frac{1}{p}-\frac{1}{\theta})\int_{\mathbb{R}^N}|\nabla
v_{n}|^{p}+(\frac{1}{p}-\frac{1}{2\theta})\int_{\mathbb{R}^N}(\lambda
A(\epsilon x)+1)|f(v_{n})|^{p},
$$
where $o_{n}(1) \to 0$ as $n \to \infty$. Recalling that $|\nabla
v_{n}|_{p}\leq 1 + |\nabla v_{n}|_{p}^{p}$,
\begin{equation}
c+o_{n}(1)\|v_{n}\|_{\epsilon,\lambda}\geq
\frac{(\theta-p)}{p\theta}\Big(|\nabla
v_{n}|_{p}-1+\int_{\mathbb{R}^N}(\lambda A(\epsilon
x)+1)|f(v_{n})|^{p}\Big)
\end{equation}
from where it follows the inequality
\begin{eqnarray*}
c_1+o_{n}(1)\|v_{n}\|_{\epsilon,\lambda}&\geq &
\frac{(\theta-p)}{p\theta} \Big(|\nabla
v_{n}|_{p}+1+\int_{\mathbb{R}^N}(\lambda A(\epsilon
x)+1)|f(v_{n})|^{p}\Big)\\
&\geq &\frac{(\theta-p)}{p\theta}\|v_{n}\|_{\epsilon,\lambda}.
\end{eqnarray*}
Thus,
$$
\limsup_{n \to \infty}\|v_{n}\|_{\epsilon,\lambda} \leq
c_1\frac{p\theta}{(\theta-p)}:=K.
$$
\hfill\rule{2mm}{2mm}

\begin{lemma} \label{L2}
Suppose that $h$ satisfies $(H_{1})-(H_{3})$. Let $(v_{n}) \subset
X_{\epsilon,\lambda}$ be a $(PS)_c$ sequence for
$I_{\epsilon,\lambda}$. Then $c \geq 0$, and if $c=0$, we have that
$v_{n} \to 0$ in $X_{\epsilon,\lambda}$.
\end{lemma}
\noindent {\bf Proof}
As in the proof of Lemma \ref{L1},
\begin{equation} \label{E1}
c+o_{n}(1)\|v_{n}\|_{\epsilon,\lambda}=I_{\epsilon,\lambda}(v_{n})-\frac{1}{\theta}I'_{\epsilon,\lambda}(v_{n})v_{n}
\geq
\Big(\frac{1}{p}-\frac{1}{\theta}\Big)\|v_{n}\|_{\epsilon,\lambda}\geq
0
\end{equation}
that is
$$
c+o_{n}(1)\|v_{n}\|_{\epsilon,\lambda}\geq 0.
$$
The boundedness of $(v_{n})$ in $X_{\epsilon,\lambda}$ gives $c \geq
0$ after passage to the limit as $n \to \infty$. If $c=0$, the inequality (\ref{E1}) gives  $v_{n} \to 0$ in $X_{
\epsilon,\lambda}$ as $n \to \infty$, finishing the proof of Lemma \ref{L2}.
\hfill\rule{2mm}{2mm}

\begin{lemma} \label{L3}
Suppose that $h$ satisfies $(H_{1})-(H_{3})$. Let $c>0$ and
$(v_{n})$ be a $(PS)_{c}$ sequence for $I_{\epsilon,\lambda}$. Then,
there exists $\delta>0$ such that
$$
\liminf_{n \to \infty}\int_{\mathbb{R}^{N}}|f(v_{n})|^{q} \geq
\delta,
$$
with $\delta$ being independent of $\lambda$ and $\epsilon$.
\end{lemma}
\noindent {\bf Proof}
From $(H_{1})-(H_{2})$, there exists a constant $C>0$ such that
\begin{equation} \label{E2}
|h(t)t|\leq \frac{1}{4}|t|^{p}+C|t|^{q} \,\,\, \forall t \in \mathbb{R}.
\end{equation}
Now, $I'_{\epsilon,\lambda}(v_{n})v_{n}=o_{n}(1)$ and  Lemma \ref{Lema f}$(6)$ give

\begin{equation} \label{E3}
\frac{1}{2}\biggl[\displaystyle\int_{\mathbb{R}^{N}}|\nabla
v_{n}|^{p} + \displaystyle\int_{\mathbb{R}^{N}}(\lambda A(\epsilon
x)+1)|f(v_{n})|^{p} \biggl]\leq \displaystyle\int_{\mathbb{R}^{N}}
h(f(v_{n})f(v_{n}).
\end{equation}
Combining (\ref{E2}) with (\ref{E3}),
\begin{equation} \label{E4}
\frac{1}{4}\biggl[\displaystyle\int_{\mathbb{R}^{N}}|\nabla
v_{n}|^{p} + \displaystyle\int_{\mathbb{R}^{N}}(\lambda A(\epsilon
x)+1)|f(v_{n})|^{p} \biggl]\leq C
\displaystyle\int_{\mathbb{R}^{N}}|f(v_{n})|^{q}.
\end{equation}
On the other hand, we have the equality
$$
\frac{1}{p}\biggl[\displaystyle\int_{\mathbb{R}^{N}}| \nabla
v_{n}|^{p} + \displaystyle\int_{\mathbb{R}^{N}}(\lambda A(\epsilon
x)+1)|f(v_{n})|^{p} \biggl]=I_{
\epsilon,\lambda}(v_{n})+\int_{\mathbb{R}^{N}}H(f(v_{n}))
$$
which combined with $(H_{3})$ and $I_{\epsilon,\lambda}(v_{n})=c+o_{n}(1)$ leads to
\begin{equation} \label{E5}
\liminf_{n \to \infty} \biggl[ \displaystyle\int_{\mathbb{R}^{N}}|\nabla
v_{n}|^{p} + \displaystyle\int_{\mathbb{R}^{N}}(\lambda A(\epsilon
x)+1)|f(v_{n})|^{p}\biggl] \geq pc>0.
\end{equation}
Now, the lemma follows from (\ref{E4}) and (\ref{E5}).
\hfill\rule{2mm}{2mm}

\begin{lemma} \label{L4}
Suppose that $h$ satisfies $(H_{1})-(H_{3})$ and $A$ satisfies
$(A_{1})-(A_{2})$. Let $d>0$ be an arbitrary number. Given any
$\epsilon>0$ and $\eta>0$, there exist $\Lambda_{\eta}>0$ and
$R_{\eta}>0$, which are independent of $\epsilon$, such that if
$(v_{n})$ is a $(PS)_{c}$ sequence for $I_{\epsilon,\lambda}$ with
$c \leq d$ and $\lambda \geq \Lambda_{\eta}$, then
$$
\limsup_{n \to \infty}\int_{\mathbb{R}^{N} \setminus
B_{R_{\eta}}(0)}|f(v_{n})|^{q}<\eta.
$$
\end{lemma}
\noindent {\bf Proof.} \, Given any $R>0$, define
$$
X(R)=\{x \in \mathbb{R}^{N}: \, |x|>R; \, A(\epsilon x)\geq K_{0}\}
$$
and
$$
Y(R)=\{x \in \mathbb{R}^{N}:\, |x|>R; \, A(\epsilon x)<K_{0}\}.
$$
Observe that
$$
\int_{X(R)}|f(v_{n})|^{p}\leq \frac{1}{\lambda
K_{0}+1}\int_{X(R)}(\lambda A(\epsilon x)+1)|f(v_{n})|^{p}.
$$
From Lemma \ref{L1}, there exists $K>0$ such that
\begin{equation} \label{E6}
\limsup_{n \to \infty}\int_{X(R)}|f(v_{n})|^{p}\leq \frac{K}{\lambda
K_{0}+1}.
\end{equation}
On the other hand, by H\"{o}lder inequality
$$
\int_{Y(R)}|f(v_{n})|^{p} \leq
\Big(\int_{Y(R)}|f(v_{n})|^{p^{*}}\Big)^{\frac{p}{p^{*}}}(\mu(Y(R)))^{\frac{p}{N}}.
$$
Using Sobolev embeddings together with Lemmas \ref{Lema f} and \ref{L1}, there exists a
constant $\widehat{K}>0$ such that
\begin{equation} \label{E7}
\limsup_{n \to \infty}\int_{Y(R)}|f(v_{n})|^{p} \leq \widehat{K}
(\mu(Y(R)))^{\frac{p}{N}},
\end{equation}
where the constant $\widehat{K}$ is uniform on $c \in [0,d]$. Since
$$
Y(R) \subset \{x \in \mathbb{R}^{N}:\, A(\epsilon x) \leq K_{0}\}
$$
it follows from $(A_{2})$
\begin{equation}\label{E8}
\lim_{R \to \infty}\mu(Y(R))=0.
\end{equation}
Using interpolation,
$$
|f(v_{n})|_{L^{q}(\mathbb{R}^{N} \setminus B_{R}(0))} \leq
|f(v_{n})|^{\alpha}_{L^{p}(\mathbb{R}^{N} \setminus
B_{R}(0))}|f(v_{n})|^{1-\alpha}_{L^{p^{*}}(\mathbb{R}^{N} \setminus
B_{R}(0))}
$$
for some $\alpha \in (0,1)$. Then, by Lemma \ref{L1}, there
exists a constant $\widetilde{K}>0$ such that
\begin{equation}\label{E9}
\limsup_{n \to \infty}\int_{\mathbb{R}^{N} \setminus
B_{R}(0)}|f(v_{n})|^{q} \leq \widetilde{K }\limsup_{n \to
\infty}\Big(\int_{\mathbb{R}^{N} \setminus
B_{R}(0)}|f(v_{n})|^{p}\Big)^{\frac{q\alpha}{p}}.
\end{equation}
Combining (\ref{E6}) with (\ref{E7}) and (\ref{E8}), given $\eta>0$,
we can fix $R=R_{\eta}$ and $\Lambda_{\eta}>0$ such that
\begin{equation}\label{E10}
\limsup_{n \to \infty}\int_{\mathbb{R}^{N} \setminus
B_{R}(0)}|f(v_{n})|^{p} \leq
\Big(\frac{\eta}{2\widetilde{K}}\Big)^{\frac{p}{q\alpha}}
\end{equation}
for all $\lambda \geq \Lambda_{\eta}$. Consequently, from (\ref{E9})
and (\ref{E10}),
$$
\limsup_{n \to \infty} \int_{\mathbb{R}^{N} \setminus
B_{R}(0)}|f(v_{n})|^{q} \leq \eta.
$$
concluding the proof of the lemma. \hfill\rule{2mm}{2mm}

\vspace{0.5 cm}

As a first consequence of the last lemma, we have the following result

\begin{corollary}\label{C1}
If $(v_{n})$ is a $(PS)_c$ sequence for $I_{\epsilon,\lambda}$ and
$\lambda$ is large enough, then its weak limit is nontrivial
provided that $c>0$.
\end{corollary}

\vspace{0.3 cm}

The next result shows that $I_{\epsilon,\lambda}$ satisfies the Palais-Smale condition for
$\lambda$ sufficiently large for $\epsilon$ arbitrary.

\vspace{0.3 cm}

\begin{proposition} \label{P1}
Suppose that $(H_{1})-(H_{3})$ and $(A_{1})-(A_{2})$ hold. Then for
any $d>0$ and $\epsilon >0$ there exists $\Lambda>0$, independent of
$\epsilon$, such that $I_{ \epsilon,\lambda}$ satisfies the
$(PS)_{c}$ condition for all $c \leq d, \lambda \geq \Lambda$ and
$\epsilon
>0$. That is, any sequence $(v_{n}) \subset X_{\epsilon,\lambda}$
satisfying
\begin{equation}
I_{\epsilon,\lambda}(v_{n}) \to c \,\,\, \mbox{and} \,\,\, I'_{
\epsilon,\lambda}(v_{n}) \to 0, \label{ps}
\end{equation}
for $c \leq d$, has a strongly convergent subsequence in $X_{
\epsilon,\lambda}$.
\end{proposition}
\noindent {\bf Proof}
Given any $d>0$ and $\epsilon>0$, take $c \leq d$ and let $(v_{n})$
be a $(PS)_{c}$ sequence for $I_{\epsilon,\lambda}$. From Lemma
\ref{L1}, there are a subsequence still denoted by $(v_{n})$ and $v
\in X_{\epsilon,\lambda}$ such that $(v_{n})$ is weakly convergent
to $v$ in $X_{\epsilon,\lambda}$. If $\widetilde{v}_{n}=v_{n}-v$,
arguing as in \cite[Lemma 3.7]{AlvesFigueiredoSevero}, it follows
that
\begin{equation} \label{E11}
I_{\epsilon,\lambda}(\widetilde{v}_{n})=I_{
\epsilon,\lambda}(v_{n})-I_{\epsilon,\lambda}(v)+o_{n}(1)
\end{equation}
and
\begin{equation} \label{E12}
I'_{\epsilon,\lambda}(\widetilde{v}_{n}) \to 0.
\end{equation}
Once that $I'_{ \epsilon,\lambda}(v)=0$, $(H_{3})$ gives

\begin{equation} \label{E13}
I_{ \epsilon,\lambda}(v)=I_{
\epsilon,\lambda}(v)-\frac{1}{\theta}I'_{ \epsilon,\lambda}(v)v \geq
\Big(\frac{1}{p}-\frac{1}{\theta}\Big)\|v\|_{\epsilon,\lambda,}\geq
0.
\end{equation}
Setting $c'=c-I_{\epsilon,\lambda}(v)$, by (\ref{E11})-(\ref{E13}),
we deduce that $c' \leq d$ and $(\widetilde{v}_{n})$ is a
$(PS)_{c'}$ sequence for $I_{ \epsilon,\lambda,}$, thus by Lemma
\ref{L2}, we have $c' \geq 0$. We claim that $c'=0$. On the
contrary, suppose that $c'>0$. From Lemma \ref{L3}, there is
$\delta>0$ such that
\begin{equation} \label{E14}
\liminf_{n \to
\infty}\int_{\mathbb{R}^{N}}|f(\widetilde{v}_{n})|^{q}>\delta.
\end{equation}
Letting $\eta=\frac{\delta}{2}$ and applying Lemma \ref{L4}, we get
$\Lambda >0$ and $R>0$ such that
\begin{equation} \label{E15}
\limsup_{n \to \infty}\int_{\mathbb{R}^{N} \setminus
B_{R}(0)}|f(\widetilde{v}_{n})|^{q} < \frac{\delta}{2}
\end{equation}
for the corresponding $(PS)_{c'}$ sequence for $I_{
\epsilon,\lambda}$ for all $\lambda \geq \Lambda$. Combining
(\ref{E14}) with (\ref{E15}) and using the fact that
$\widetilde{v}_{n} \rightharpoonup 0$ in $X_{ \epsilon,\lambda}$, we
derive
$$
\delta \leq \liminf_{n \to
\infty}\int_{\mathbb{R}^{N}}|f(\widetilde{v}_{n})|^{q} \leq \limsup_{n \to
\infty}\int_{\mathbb{R}^{N} \setminus
B_{R}(0)}|f(\widetilde{v}_{n})|^{q}\leq \frac{\delta}{2}
$$
which is impossible, then $c'=0$. Thereby,  by Lemma \ref{L2}, $\widetilde{v}_{n} \to 0$ in $X_{\epsilon,\lambda}$, that is,
$v_{n} \to v$ in $X_{\epsilon,\lambda}$ and the proof of Proposition
\ref{P1} is complete.
\hfill\rule{2mm}{2mm}

\vspace{0,2cm}

In closing this section, we proceed with the study of
$(PS)_{c,\infty}$ sequences, that is, sequences $(v_{n})$ in
$X_{\epsilon,\lambda}$ verifying:
$$
\begin{array}{l}
i)\,\,\,\,\,\, \lambda_{n} \to \infty \\
\mbox{}\\
ii)\,\,\,\,\,\, (I_{\epsilon,\lambda_{n}}(v_{n})) \,\, \mbox{is bounded} \\
\mbox{}\\
iii)\,\,\,\,\,\, \|I'_{\epsilon,\lambda_{n}}(v_{n})\|^{*}_{\epsilon,\lambda_{n}} \to 0\\
\end{array}
$$
where $\|\,\,\,\|^{*}_{\epsilon,\lambda_{n}}$ is defined by
$$
\|\varphi\|^{*}_{\epsilon,\lambda_{n}}= \sup \{|\varphi(u)|;\  u \in
X_{\epsilon,\lambda_{n}},\,  \|u\|_{\epsilon,\lambda}\leq 1\}
 \,\,\,\, \mbox{for} \,\,\,\, \varphi \in X^{*}_{\epsilon,\lambda_{n}}.
$$

\begin{proposition}\label{P2}
Suppose that $(H_{1})-(H_{3})$ and $(A_{1})-(A_{2})$ hold. Assume
that $(v_{n}) \subset W^{1,p}(\mathbb{R}^{N})$ is a
$(PS)_{c,\infty}$ sequence. Then for each $\epsilon >0$ fixed, there
exists a subsequence still denoted by $(v_{n})$ and $v_{\epsilon}
\in W^{1,p}(\mathbb{R}^{N})$ such that \vspace{0.2 cm}
\begin{itemize}
\item[i)]\ $v_{n} \to v_{\epsilon}$ in $W^{1,p}(\mathbb{R}^N)$. Moreover, $v_{\epsilon}=0$ on $\Omega_{\epsilon}^{c}$ and
$v_{\epsilon} \in W^{1,p}(\mathbb{R}^N)\bigcap
L^{\infty}_{loc}(\mathbb{R}^N)$ is a solution of
$$
\left\{
\begin{array}{l}
-\Delta_{p}v+|f(v)|^{p-2}f(v)f'(v)=h(f(v))f'(v), \,\,\, \mbox{in} \,\,\, \Omega_{\epsilon} \\
\mbox{}\\
v>0 \,\, \mbox{in} \,\, \Omega  \,\,\ \mbox{and} \,\,\, v=0 \,\,\,
\mbox{on} \,\,
\partial \Omega_{\epsilon}
\end{array}
\right.
$$
where $\Omega_{\epsilon}=\frac{\Omega}{\epsilon}$.

\item[ii)]\ $\displaystyle
\lambda_{n}\int_{\mathbb{R}^{N}}|f(v_{n})|^{p} \to 0$.

\item[iii)]\ $\displaystyle \|v_{n}-v\|_{\epsilon,\lambda_{n}}
\to 0$.
\end{itemize}
\end{proposition}

\vspace{0.2 cm}

\noindent \textbf{Proof.} As in the proof of Lemma \ref{L1}, the
sequence $(\|v_{n}\|_{\epsilon,\lambda_{n}})$ is bounded in
$\mathbb{R}$. Thus, we can extract a subsequence $v_{n}
\rightharpoonup v_{\epsilon}$ weakly in $X_{\epsilon,\lambda}$. For each
$m \in \mathbb{N}$, we define the set
$$
C_{m}=\Big\{x \in \mathbb{R}^{N}:\, A_{\epsilon}(x)\geq
\frac{1}{m}\Big\}, \,\,\, \mbox{where} \,\,\,
A_{\epsilon}(x)=A(\epsilon x)
$$
which satisfies
$$
\int_{C_{m}}|f(v_{n})|^{p} \leq
m\int_{C_{m}}A_{\epsilon}(x)|f(v_{n})|^{p}\leq
\frac{m}{\lambda_{n}}\int_{\mathbb{R}^{N}}(1+\lambda_{n}A_{\epsilon}(x))|f(v_{n})|^{p}.
$$
Thus, from Lemma \ref{L1},
$$
\int_{C_{m}}|f(v_{n})|^{p} \leq \frac{mK}{\lambda_{n}} \,\,\,
\mbox{for} \,\,\, n \in \mathbb{N},
$$
for some constant $K>0$. Hence by Fatou's Lemma,
$$
\int_{C_{m}}|f(v_{\epsilon})|^{p}=0
$$
after to passage to the limit as $n \to \infty$. Thus $f(v_{\epsilon})=0$ almost everywhere in $C_{m}$. Once that $f(t)=0$ if, and only if $t=0$, it follows that $v_{\epsilon}=0$ almost everywhere in $C_{m}$. Observing that
$$
\mathbb{R}^{N} \setminus
A_{\epsilon}^{-1}(0)=\cup_{m=1}^{\infty}C_{m},
$$
we deduce that $v_{\epsilon}=0$ almost everywhere in $\mathbb{R}^{N}
\setminus A_{\epsilon}^{-1}(0)$.
 Now, recalling that $A_{\epsilon}^{-1}(0)=\overline{\Omega}_{\epsilon} \cup D_{\epsilon}$ and
 $\mu(D_{\epsilon})= \mu(\frac{1}{\epsilon}D)=0$, it follows that $v_{\epsilon}=0$ almost everywhere
 in
 $\mathbb{R}^{N} \setminus \overline{\Omega}_{\epsilon}$. As $\partial \Omega_{\epsilon}$ is a smooth set, let us conclude that
 $v_{\epsilon} \in W^{1,p}_{0}(\Omega_{\epsilon})$.

Arguing as in Lemma \ref{L4}, we can assert that given any $\eta>0$
there exists $R>0$ such that
\begin{equation}
\limsup_{n \to \infty}\int_{\mathbb{R}^{N} \setminus
B_{R}(0)}|f(v_{n})|^{p} < \eta \label{lim1}
\end{equation}
 and
\begin{equation}
\limsup_{n \to \infty}\int_{\mathbb{R}^{N} \setminus
B_{R}(0)}|f(v_{n})|^{q} < \eta. \label{lim2}
\end{equation}
From $(H_{1})-(H_{2})$, for each $\tau >0$ there exists $C_{\tau}>0$
such that
$$
|h(s)|\leq \tau |s|^{p-1}+C_{\tau}|s|^{q-1} \,\,\, \mbox{for all}
\,\,\, s \in \mathbb{R}.
$$
This inequality combined with Sobolev's embeddings and the
limits ({\ref{lim1}) and (\ref{lim2}) yields  there is  a
subsequence, still denoted by $(v_{n})$, such that
\begin{equation} \label{E16}
\lim_{n \to \infty}\int_{\mathbb{R}^{N}}h(f(v_{n}))f'(v_{n})v_{n}=\int_{\mathbb{R}^{N}}h(f(v_{\epsilon}))f'(v_{\epsilon})v_{\epsilon},
\end{equation}
and
\begin{equation} \label{E17}
\lim_{n \to \infty}\int_{\mathbb{R}^{N}}h(f(v_{n}))f'(v_{\epsilon})v_{\epsilon}=\int_{\mathbb{R}^{N}}h(f(v_{\epsilon}))f'(v_{\epsilon})v_{\epsilon}
\end{equation}
In the sequel, define
$$
\begin{array}{l}
P_{n}=\displaystyle \int_{\mathbb{R}^{N}}\left\langle |\nabla v_{n}|^{p-2}\nabla v_{n}- |\nabla v_{\epsilon}|^{p-2}\nabla v_{\epsilon}, \nabla v_{n} -
\nabla v_{\epsilon}  \right\rangle  + \\
\mbox{}\\
\hspace{1 cm} +\displaystyle \int_{\mathbb{R}^{N}}( |f(v_{n})|^{p-2}f(v_{n})f'(v_{n})-
|f(v_{\epsilon})|^{p-2}f(v_{\epsilon})f'(v_{\epsilon}))(v_{n}-v_{\epsilon})
\end{array}
$$
and observe that
$$
\begin{array}{l}
P_{n} \leq \displaystyle \int_{\mathbb{R}^{N}}\left\langle |\nabla v_{n}|^{p-2}\nabla v_{n}-
|\nabla v_{\epsilon}|^{p-2}\nabla v_{\epsilon}, \nabla v_{n} - \nabla v_{\epsilon}  \right\rangle + \\
\mbox{}\\
\displaystyle \int_{\mathbb{R}^{N}}(\lambda_{n}A_{\epsilon}(x)+1) ( |f(v_{n})|^{p-2}f(v_{n})f'(v_{n})-
|f(v_{\epsilon})|^{p-2}f(v_{\epsilon})f'(v_{\epsilon}))(v_{n}-v_{\epsilon})=  \\
\mbox{}\\
\hspace{1 cm} I'_{\epsilon,\lambda_{n},}(v_{n})v_{n}-
I'_{\epsilon,\lambda_{n}}(v_{n})v_{\epsilon}+ \displaystyle
\int_{\mathbb{R}^{N}}h(f(v_{n}))f(|v_{n})|^{p-2}f(v_{n})f'(v_{n})v_{n}\\
\hspace{1 cm} -
\displaystyle\int_{\mathbb{R}^{N}}h(f(v_{n}))f(|v_{n})|^{p-2}f(v_{n})f'(v_{n})v_{\epsilon}+o_{n}(1).
\end{array}
$$
Thus, by (\ref{E16}) and (\ref{E17}) it follows that
$P_{n}=o_{n}(1),$ that is,
\begin{equation} \label{TT1}
\int_{\mathbb{R}^{N}}\left\langle |\nabla v_n|^{p-2}\nabla v_n - |\nabla v_\epsilon|^{p-2}\nabla v_\epsilon, \nabla v_n - \nabla v_\epsilon \right\rangle = o_n(1)
\end{equation}
and
\begin{equation} \label{E172}
\int_{\mathbb{R}^{N}}(|f(v_n)|^{p-2}f(v_n)f'(v_n)-|f(v_\epsilon)|^{p-2}f(v_\epsilon)f'(v_\epsilon))(v_n-v_\epsilon) = o_n(1).
\end{equation}
The limits (\ref{E172}) and (\ref{TT1}) combined with Corollary \ref{monotonia} give
\begin{equation} \label{E19}
v_{n} \to v_{\epsilon} \,\,\, \mbox{strongly in} \,\,\,
W^{1,p}(\mathbb{R}^{N}).
\end{equation}
Now, using the fact that $(\lambda_{n}A_{\epsilon}(x)+1)v_{\epsilon}(x)=v_{\epsilon}(x)$ a.e in $\mathbb{R}^{N}$ and that for each $\phi \in C^{\infty}_{0}(\Omega_{\epsilon})$, $I'_{\epsilon,\lambda_{n}}(v_{n})\phi=o_{n}(1)$, we have that
$$
\int_{\mathbb{R}^{N}}(|\nabla v_{n}|^{p-2}\nabla v_{n} \nabla  \phi
+ |f(v_{n})|^{p-2}f(v_{n})f'(v_{n})
\phi)=\int_{\mathbb{R}^{N}}h(f(v_{n}))f'(v_{n})\phi +o_{n}(1).
$$
This together with (\ref{E19}) yields
$$
\int_{\Omega_{\epsilon}}(|\nabla v_{\epsilon}|^{p-2}\nabla
v_{\epsilon} \nabla  \phi +
|f(v_{\epsilon})|^{p-2}f(v_{\epsilon})f'(v_{\epsilon})
\phi)=\int_{\Omega_{\epsilon}}h(f(v_{\epsilon}))f'(v_{\epsilon})\phi
$$
and hence
\begin{equation} \label{E20}
\int_{\Omega_{\epsilon}}(|\nabla v_{\epsilon}|^{p-2}\nabla
v_{\epsilon} \nabla  w +
 |f(v_{\epsilon})|^{p-2}f(v_{\epsilon})f'(v_{\epsilon})w)=\int_{\Omega_{\epsilon}}h(f(v_{\epsilon}))f'(v_{\epsilon})w,
\end{equation}
for all $w \in W_{0}^{1,p}(\Omega_{\epsilon})$. Arguing as \cite[Proposition 3.6]{AlvesFigueiredoSevero},  we can
prove that $v_{\epsilon} \in L^{\infty}(\mathbb{R}^{N})$. Thereby, $v_{\epsilon}$ is a solution of
$$
\left\{
\begin{array}{l}
-\Delta_{p}{v}+|f(v)|^{p-2}f(v)f'(v)=h(f(v))f'(v), \,\,\, \mbox{in} \,\,\, \Omega_{\epsilon} \\
\mbox{}\\
v>0 \,\, \mbox{in} \,\, \Omega  \,\,\ \mbox{and} \,\,\, u=0 \,\,\,
\mbox{on} \,\,\partial \Omega_{\epsilon}
\end{array}
\right.
$$
and the proof of $i)$ is complete.

To deduce $ii)$, we start observing that
\begin{eqnarray*}
&&\displaystyle\int_{\mathbb{R}^{N}}|\nabla v_{n}|^{p} +
\displaystyle\int_{\mathbb{R}^{N}}|f(v_{n})|^{p-2}f(v_{n})f'(v_{n})
v_{n}\\
&+&
\lambda_{n}\displaystyle\int_{\mathbb{R}^{N}}A_{\epsilon}|f(v_{n})|^{p-2}f(v_{n})f'(v_{n})
v_{n}=\displaystyle\int_{\mathbb{R}^{N}}h(f(v_{n}))f'(v_{n})v_{n}
+o_{n}(1).
\end{eqnarray*}

The last equality, combined with (\ref{E19}) and (\ref{E20}) leads to
$$
\lim_{n \to \infty}
\lambda_{n}\displaystyle\int_{\mathbb{R}^{N}}A_{\epsilon}(x)|f(v_{n})|^{p-2}f(v_{n})f'(v_{n})
v_{n}=0.
$$
This limit together with Lemma \ref{Lema f}$(6)$ implies that
$$
\displaystyle\lim_{n \to \infty}
\lambda_{n}\int_{\mathbb{R}^{N}}A_{\epsilon}(x)|f(v_{n})|^{p}=0,
$$
proving ii).

For to prove iii), we observe that
$$
\lambda_n \int_{\mathbb{R}^{N}}A_\epsilon(x)|f(v_n -v_{\epsilon})|^{p}=\lambda_n \int_{\mathbb{R}^{N} \setminus \Omega }A_\epsilon(x)|f(v_n)|^{p} \leq \lambda_n \int_{\mathbb{R}^{N} }A_\epsilon(x)|f(v_n)|^{p}
$$
because $v_{\epsilon}=0$ in $\Omega$. Hence,
$$
\lambda_n \int_{\mathbb{R}^{N}}A_\epsilon(x)|f(v_n -v_{\epsilon})|^{p} \to 0 \,\,\, \mbox{as} \,\,\, n \to +\infty
$$
which together with Corollary \ref{resultadoClaudianor} leads to
\begin{equation} \label{Z1}
\int_{\mathbb{R}^{N}}(1+\lambda_n A_\epsilon(x))|f(v_n -v_{\epsilon})|^{p} \to 0 \,\,\, \mbox{as} \,\,\, n \to +\infty.
\end{equation}
Recalling that
$$
\|v_{n}-v_{\epsilon}\|_{\epsilon,\lambda_{n}} \leq |\nabla v_{n}- \nabla v_{\epsilon}|_{p} + \int_{\mathbb{R}^{N}}(1+\lambda_n A_\epsilon(x))|f(v_n -v_{\epsilon})|^{p},
$$
it follows from (\ref{E19}) and (\ref{Z1}),
$$
\lim_{n \to \infty}\|v_{n}-v_{\epsilon}\|_{\epsilon,\lambda_{n}}=0,
$$
which proves $iii)$, and the proof of Proposition \ref{P2} is
complete. $\hfill \rule{2mm}{2mm}$

\begin{corollary} \label{T2}
Suppose that $(A_{1})-(A_{2})$ and $(H_{1})-(H_{4})$ hold. Then for
each $\epsilon >0$ and a sequence $(v_{n})$ of solutions of
$(P_{\epsilon,\lambda_{n}})$ with $\lambda_{n} \to \infty$ and
$\limsup_{n\to \infty}I_{\epsilon,\lambda_n}(v_{n})< \infty$, there
exists a subsequence that converges strongly in
$W^{1,p}(\mathbb{R}^{N})$ to a solution of the problem
$$
\left\{
\begin{array}{l}
-\Delta_{p}{v}+|f(v)|^{p-2}f(v)f'(v)=h(f(v))f'(v), \,\,\, \mbox{in} \,\,\, \Omega_\epsilon \\
\mbox{}\\
v>0 \,\, \mbox{in} \,\, \Omega_\epsilon  \,\,\ \mbox{and} \,\,\, v=0
\,\,\,
\partial \Omega_\epsilon.
\end{array}
\right.
$$
\end{corollary}

\noindent \textbf{Proof.} By assumptions, there exist $c\in
\mathbb{R}$ and a subsequence of $(v_n)$, still denoted by $(v_n)$,
such that $(v_n)$ is a $(PS)_{c,\infty}$ sequence. The rest of the
proof follows from Proposition \ref{P2}. $\hfill \rule{2mm}{2mm}$

%%%%%%%%%%%%%%%%%%%%%%%%%%%%%%%%%%%%%% EBED %%%%%%%%%%%%%%%%%%%%%%%%%%%%%%%%%%%%%

\section{Behavior of minimax levels}

This section is devoted to the study of the behavior of the minimax
levels with respect to parameter $\lambda$ and $\epsilon$. For this
purpose, we introduce some notations. In the next,
$\mathcal{M}_{\epsilon,\lambda}$ denotes the Nehari manifold
associated to $I_{\epsilon,\lambda}$, that is,
$$
\mathcal{M}_{\epsilon,\lambda}=\Big\{v \in X_{\epsilon,\lambda}:\,
v\not=0 \,\,\, \mbox{and} \,\,\, I'_{\epsilon,\lambda}(v)v=0\Big\}
$$
and
$$
c_{\epsilon,\lambda}=\inf_{v \in
\mathcal{M}_{\epsilon,\lambda}}I_{\epsilon,\lambda}(v).
$$
From $(H_1)-(H_4)$, as proved in \cite[Lemma
3.3]{AlvesFigueiredoSevero}, the number $c_{\epsilon,\lambda}$ is
the mountain pass minimax level associated with
$I_{\epsilon,\lambda}$.

On account of the proof of Proposition \ref{P2}, when $\lambda$ is
large, the following problem can be seen as a limit problem of
$(D_{\epsilon,\lambda,})$ for each $\epsilon>0$:
$$
\left\{
\begin{array}{l}
-\Delta_{p}{v}+|f(v)|^{p-2}f(v)f'(v)=h(f(v))f'(v), \,\,\, \mbox{in} \,\,\, \Omega_{\epsilon} \\
\mbox{}\\
v>0 \,\, \mbox{in} \,\, \Omega  \,\,\ \mbox{and} \,\,\, v=0 \,\,\,
\mbox{on} \,\, \partial \Omega_{\epsilon}
\end{array}
\right. \leqno{(D_{\epsilon})}
$$
whose corresponding functional is given by
$$
E_{\epsilon}(v)=\frac{1}{p}\int_{\Omega_{\epsilon}}(|\nabla
v|^{p}+|f(v)|^{p})-\int_{\Omega_{\epsilon}}H(f(v))
$$
for every $v \in W_{0}^{1,p}(\Omega_{\epsilon})$. Here and
subsequently, $\mathcal{M}_{\epsilon}$ denotes the Nehari manifold
associated to $E_{\epsilon}$ and
$$
c(\epsilon,\Omega)=\displaystyle \inf_{v \in
\mathcal{M}_{\epsilon}}E_{\epsilon}(v)
$$
stands for the mountain pass minimax associated with $E_{\epsilon}$.
Since $0 \in \Omega$, there is $r>0$ such that $B_{r}=B_{r}(0)
\subset \Omega$ and
$B_{\frac{r}{\epsilon}}=B_{\frac{r}{\epsilon}}(0) \subset
\Omega_{\epsilon}$. We will denote by
$E_{\epsilon,B_{r}}:W_{0}^{1,p}(B_{\frac{r}{\epsilon}}(0)) \to
\mathbb{R}$ the functional
$$
E_{\epsilon,B_{r}}(v)=\frac{1}{p}\int_{B_{\frac{r}{\epsilon}}}(|\nabla
v|^{p}+|f(v)|^{p})-\int_{B_{\frac{r}{\epsilon}}}H(f(v)).
$$
Furthermore, we write $\mathcal{M}_{\epsilon,B_{r}}$ the Nehari
manifold associated to $E_{\epsilon,B_{r}}$ and
$$
c(\epsilon,B_{r})=\displaystyle \inf_{v \in
\mathcal{M}_{\epsilon,B_{r}}}E_{\epsilon,B_r}(v).
$$
Once that $B_{\frac{r}{\epsilon}} \subset
\Omega_{\epsilon}$, we have $ c(\epsilon,\Omega) \leq
c(\epsilon,B_{r})$ for every $\epsilon >0$.

Here it is important the number $c_{\infty}$, which denotes the
mountain minimax value associated to
$$
I_{\infty}(v)=\frac{1}{p}\int_{\mathbb{R}^{N}}(|\nabla
v|^{p}+|f(v)|^{p})-\int_{\mathbb{R}^{N}}H(f(v)) \,\,\, \mbox{for
all} \,\, v \in  W^{1,p}_{0}(\mathbb{R}^{N}),
$$
whose existence is guaranteed by \cite[Lemma 3.1]{AlvesFigueiredoSevero1}. Since $I_{\epsilon,\lambda}(tv) \geq I_{\infty}(t v)$ for all $t >0$
and $v \in W^{1,p}(\mathbb{R}^{N})$,
$$
c_{\epsilon,\lambda} \geq c_{\infty}.
$$

\vspace{0.5 cm}

\begin{proposition} \label{P3}
Suppose $(H_{1})-(H_{4})$ and $(A_{1})-(A_{2})$ hold. Let $\epsilon
>0$ be an arbitrary number. Then,
$$
\lim_{\lambda \to \infty}c_{\epsilon,\lambda}=c(\epsilon,\Omega).
$$
\end{proposition}

\noindent \textbf{Proof.} By Proposition \ref{P1} and Mountain Pass
Theorem, we can assume that there are two sequences,  $\lambda_{n}
\to \infty$ and  $(v_{n}) \subset X_{\epsilon,\lambda_{n}}$,   such
that
$$
I_{\epsilon,\lambda_{n}}(v_{n})=c_{\epsilon,\lambda_{n}}>0\,\,\,
\mbox{and} \,\,\, I'_{\epsilon,\lambda_{n}}(v_{n})=0.
$$
From definitions of $c_{\epsilon,\lambda_{n}}$ and $c(\epsilon,
\Omega)$,
$$
c_{\epsilon,\lambda_{n}} \leq c(\epsilon, \Omega) \,\,\, \mbox{for
all} \,\,\, n \in \mathbb{N}
$$
which implies
$$
0\leq I_{\epsilon,\lambda_{n}}(v_{n}) \leq c(\epsilon, \Omega)
\,\,\, \mbox{and} \,\,\, I'_{\epsilon,\lambda_{n}}(v_{n})=0.
$$
Thus, for some subsequence $(v_{n_{j}})$, there exists $c \in
[0,c(\epsilon, \Omega)]$ such that
$$
I_{\epsilon,\lambda_{n_{j}}}(v_{n_{j}})=c_{\epsilon,\lambda_{n_{j}}}
\to c \,\,\, \mbox{and} \,\,\,
I_{\epsilon,\lambda_{n_{j}}}'(v_{n_{j}}) \to 0
$$
showing that $(v_{n_{j}})$ is a $(PS)_{c,\infty}$, and so,
$$
\displaystyle\int_{\mathbb{R}^{N}}|\nabla v_{n}|^{p}
+\displaystyle\int_{\mathbb{R}^{N}}(\lambda_{n}A_{\epsilon}(x)+1)|f(v_{n})|^{p}
\geq p c_{\epsilon, \lambda n} \geq p c_{\infty} >0 \,\,\, \forall n \in \mathbb{N}.
$$
By Proposition \ref{P2},
$$
\lambda_{n}\int_{\mathbb{R}^{N}}A_{\epsilon}(x)|f(v_{n})|^{p} \to 0 \,\,\, \mbox{as} \,\,\, n \to +\infty
$$
then,
\begin{equation}\label{bom}
\displaystyle\int_{\mathbb{R}^{N}}|\nabla v_{n}|^{p}
+\displaystyle\int_{\mathbb{R}^{N}}|f(v_{n})|^{p}
 \geq p c_{\infty} >0 +o_n(1) \,\,\, \forall n \in \mathbb{N},
\end{equation}
implying that any subsequence of $(v_{n})$ does not converge to zero in
$W^{1,p}(\mathbb{R}^{N})$. From Proposition \ref{P2}, there exist a subsequence
$(v_{n_{j_{k}}})$ and $v \in W^{1,p}(\mathbb{R}^{N})$ such that
\begin{equation} \label{E24}
v_{n_{j_{k}}} \to v \,\,\, \mbox{strongly in} \,\,\,
W^{1,p}(\mathbb{R}^{N}) \,\,\, \mbox{and} \,\,\, v=0 \,\,\,
\mbox{in} \,\,\, \mathbb{R}^{N} \setminus \Omega_{\epsilon}.
\end{equation}
From (\ref{bom}) and (\ref{E24}), $v \not=0$ in
$W^{1,p}_{0}(\Omega_{\epsilon})$ and $v$ satisfies
$$
\left\{
\begin{array}{l}
-\Delta_{p}{u}+|f(v)|^{p-2}f(v)f'(v)=h(f(v))f'(v), \,\,\, \mbox{in} \,\,\, \Omega_{\epsilon} \\
\mbox{}\\
v>0 \,\, \mbox{in} \,\, \Omega  \,\,\ \mbox{and} \,\,\, v=0 \,\,\,
\mbox{on} \,\,  \partial \Omega_{\epsilon},
\end{array}
\right.
$$
from where it follows that
\begin{equation} \label{E25}
E_{\epsilon}(v) \geq c(\epsilon,\Omega).
\end{equation}
On the other hand,
\begin{equation}\label{E26}
E_{\epsilon}(v)=\lim_{k\to
\infty}I_{\lambda_{n_{j_{k}}},\epsilon}(v_{n_{j_{k}}})=\lim_{k\to
\infty}c_{\epsilon,\lambda_{n_{j_{k}}}}=c\leq c(\epsilon,\Omega).
\end{equation}
Therefore, (\ref{E25}) and (\ref{E26}) give
$$
\lim_{k\to
\infty}c_{\epsilon,\lambda_{n_{j_{k}}}}=c(\epsilon,\Omega).
$$
As a result,  $c_{\epsilon,\lambda} \to c(\epsilon,\Omega)$
as $\lambda \to \infty$, and the lemma follows. $\hfill\rule{2mm}{2mm}$

\vspace{0,2cm}

\begin{corollary} \label{T3}
Suppose that $(A_{1})-(A_{2})$ and $(H_{1})-(H_{4})$ hold. Then for
each $\epsilon >0$ and a sequence $(v_{n})$ of least energy
solutions of $(D_{\epsilon,\lambda_{n}})$ with $\lambda_{n} \to
\infty$ and $\limsup_{n\to \infty}I_{\epsilon,\lambda_n}(v_{n})<
\infty$, there exists a subsequence that converges strongly in
$W^{1,p}(\mathbb{R}^{N})$ to a least energy solution of the problem
$$
\left\{
\begin{array}{l}
-\Delta_{p}{u}+|f(v)|^{p-2}f(v)f'(v)=h(f(v))f'(v), \,\,\, \mbox{in} \,\,\, \Omega_\epsilon \\
\mbox{}\\
u>0 \,\, \mbox{in} \,\, \Omega_\epsilon  \,\,\ \mbox{and} \,\,\, u=0
\,\,\,
\partial \Omega_\epsilon.
\end{array}
\right.
$$
\end{corollary}

\noindent \textbf{Proof.} The proof is a consequence of Propositions \ref{P2}
and \ref{P3}. $\hfill \rule{2mm}{2mm}$

\vspace{0.5 cm}

Hereafter, $r>0$ denotes a number such that $B_{r}(0)
\subset \Omega$ and the sets
$$
\Omega_{+}=\{x \in \mathbb{R}^{N}:\, d(x,\overline{\Omega})\leq r\}
$$
and
$$
\Omega_{-}=\{x \in \mathbb{R}^{N}:\, d(x,\partial \Omega) \geq r\}
$$
are homotopically equivalent to $\Omega$. The existence of this $r$
is given by condition $(A_{1})$. For each $v \in
W^{1,p}(\mathbb{R}^{N})$ whose positive part $v_{+}=\max\{v,0\}$ is
different from zero and has a compact support, we consider the
center mass of $v$
$$
\beta(v)= \frac{\displaystyle \int_{\mathbb{R}^{N}}x
v_{+}^{p}}{\displaystyle \int_{\mathbb{R}^{N}}v_{+}^{p}}.
$$

Consider $R>0$ such that $\Omega \subset B_{R}(0)$, thus
$\Omega_{\epsilon} \subset B_{\frac{R}{\epsilon}}(0)$, and define
the auxiliary function
$$
\xi_{\epsilon}(t)= \left\{
\begin{array}{l}
1, \,\, 0\leq t \leq \frac{R}{\epsilon}\\
\mbox{}\\
\frac{R}{\epsilon t},\,\,\,\, \frac{R}{\epsilon}\leq t.
\end{array}
\right.
$$
For $v \in W^{1,p}(\mathbb{R}^{N}),v_{+}\not=0$, define
$$
\beta_{\epsilon}(v)=\frac{\displaystyle
\int_{\mathbb{R}^{N}}x\xi_{\epsilon}(|x|)v_{+}^{p}}{\displaystyle
\int_{\mathbb{R}^{N}}v_{+}^{p}}.
$$

Now for each $y \in \mathbb{R}^{N}$ and $R>2{\rm{diam}}(\Omega)$
fix
$$
A_{\frac{R}{\epsilon},\frac{r}{\epsilon},y}=\Big\{x \in
\mathbb{R}^{N}:\, \frac{r}{\epsilon}\leq |x-y|\leq
\frac{R}{\epsilon}\Big\}.
$$
We observe that if $y \notin \frac{1}{\epsilon}\Omega_{+}$ then
$\overline{\Omega_{\epsilon}} \cap
B_{\frac{r}{\epsilon}}(y)=\emptyset$. As a consequence
\begin{equation} \label{E27}
\overline{\Omega_{\epsilon}} \subset
A_{\frac{R}{\epsilon},\frac{r}{\epsilon},y}
\end{equation}
for every $y \notin \frac{1}{\epsilon}\Omega_{+}$. Moreover, for $y
\in \mathbb{R}^{N}$, $\alpha(R,r,\epsilon,y)$ denotes the number
$$
\alpha(R,r,\epsilon,y)=\inf \Big\{\widehat{J}_{\epsilon,y}(v):\,
\beta(v)= y \,\,\, \mbox{and} \,\,\, v \in
\widehat{\mathcal{M}}_{\epsilon,y}\Big\}
$$
where
$$
\widehat{J}_{\epsilon,y}(v)=\frac{1}{p}\int_{A_{\frac{R}{\epsilon},\frac{r}{\epsilon},y}}(|\nabla
v|^{p}+|f(v)|^{p})-\int_{A_{\frac{R}{\epsilon},\frac{r}{\epsilon},y}}H(f(v))
$$
and
$$
\widehat{\mathcal{M}}_{\epsilon,y}=\Big\{v \in
W^{1,p}_{0}(A_{\frac{R}{\epsilon},\frac{r}{\epsilon},y}):v\not=0
\,\,\, \mbox{and} \,\,\, \widehat{J}_{\epsilon,y}'(v)v=0\Big\}.
$$
From now on, we will write $\alpha(R,r,\epsilon,0)$ as
$\alpha(R,r,\epsilon)$, $\widehat{J}_{\epsilon,0}$ as
$\widehat{J}_{\epsilon}$ and  $\widehat{\mathcal{M}}_{\epsilon,0}$
as $\widehat{\mathcal{M}}_{\epsilon}$.

\begin{lemma} \label{L6}
Assume that $(H_{1})-(H_{4})$ hold. Then, there exists $\epsilon^{*}>0$
such that
$$
c(\epsilon,\Omega) < \alpha(R,r,\epsilon)
$$
for every $\epsilon \in (0, \epsilon^{*})$.
\end{lemma}
\noindent \textbf{Proof.} Invoking \cite[Proposition
4.1]{AlvesFigueiredoSevero1}, we assert that
$$
\lim_{\epsilon \to 0}\alpha(R,r,\epsilon)>c_{\infty}.
$$
Thus, there exists $\epsilon_{1}>0$ such that
\begin{equation} \label{E28}
\alpha(R,r,\epsilon) > c_{\infty}+\delta
\end{equation}
for all $0< \epsilon < \epsilon_{1}$, for some $\delta>0$. On the
other hand, arguing as in \cite[Proposition
4.2]{AlvesFigueiredoSevero1},
$$
\lim_{\epsilon \to 0}c(\epsilon, B_{r})=c_{\infty}.
$$

Therefore, there exists $\epsilon_{2}>0$ such that
\begin{equation} \label{E29}
c(\epsilon, B_{r}) < c_{\infty} +\frac{\delta}{2} \,\,\, \mbox{for
all} \,\,\, 0< \epsilon < \epsilon_{2}.
\end{equation}
For $\epsilon^{*}=\min\{\epsilon_{1},\epsilon_{2}\}$,  (\ref{E28})
and (\ref{E29}) lead to
$$
c(\epsilon, B_{r}) < \alpha(R,r,\epsilon)
$$
for every $\epsilon \in (0, \epsilon^{*})$. Now, the lemma follows of the inequality \linebreak $c(\epsilon, \Omega) \leq c(\epsilon, B_{r})$. $\hfill \rule{2mm}{2mm}$

\vspace{0.5 cm}

To conclude this section,  we establish a result about the center of
mass of the function in the Nehari manifold
$\mathcal{M}_{\epsilon,\lambda}$.

\begin{lemma} \label{L7}
Suppose $(H_{1})-(H_{4})$ and $(A_{1})-(A_{2})$ hold. Let
$\epsilon^{*}>0$ given by Lemma \ref{L6}. Then for any $\epsilon \in
(0, \epsilon^{*})$, there exists $\lambda^{*}>0$ which depends on
$\epsilon$ such that
$$
\beta_{\epsilon}(v) \in \frac{1}{\epsilon}\Omega_{+}
$$
for all $\lambda > \lambda^{*},0< \epsilon < \epsilon^{*}$ and $v
\in \mathcal{M}_{\epsilon,\lambda}$ with $I_{\epsilon,\lambda}(v)
\leq c(\epsilon,B_{r})$.
\end{lemma}
\noindent \textbf{Proof.} Suppose by contradiction that there exists
a sequence $(\lambda_{n})$ with $\lambda_{n} \to \infty$ such that
$$
v_{n} \in \mathcal{M}_{\epsilon,\lambda_{n}}, \,\,
I_{\epsilon,\lambda_{n}}(v_{n}) \leq c(\epsilon, B_{r})
$$
and
\begin{equation}\label{E30}
\beta_{\epsilon}(v_{n}) \notin \frac{1}{\epsilon}\Omega_{+}.
\end{equation}
Repeating the same arguments used in the proofs of Lemma \ref{L4}
and Proposition \ref{P2},
$(\|v_n\|_{\epsilon,\lambda_n})$ is a bounded sequence in
$\mathbb{R}$ and there exists $v \in W^{1,p}(\mathbb{R}^N)$ such
that $v_n \rightharpoonup v$ weakly in $W^{1,p}(\mathbb{R}^N)$, $v =
0$ in $\mathbb{R}^N\setminus \Omega_\epsilon$ and for each $\eta>0$
there exists $R>0$ such that
$$
\limsup_{n \to \infty}\int_{\mathbb{R}^{N} \setminus
B_{R}(0)}|f(v_{n})|^{p}<\eta.
$$
This fact implies that
\[
f(v_n) \to f(v) \,\,\, \mbox{strongly in} \,\,\, L^p(\mathbb{R}^N).
\]
Hence by interpolation,
\[
f(v_n) \to f(v) \,\,\, \mbox{strongly in} \,\,\,
L^t(\mathbb{R}^N)\,\,\, \mbox{for all} \,\,\, t\in [p,p^*).
\]
On the other hand, since $v_n \in
\mathcal{M}_{\epsilon,\lambda_{n}}$,  from  (\ref{bom}),
$$
0<pc_{\infty} \leq
\displaystyle\int_{\mathbb{R}^{N}}h(f(v_n))f'(v_n)v_n, \,\,\,
\mbox{for all} \,\, n \in \mathbb{N},
$$
from where it follows that
$$
0< p c_{\infty} \leq \int_{\mathbb{R}^N}h(f(v))f'(v)v,
$$
which yields
\begin{equation} \label{E31}
v \not=0, \, E'_{\epsilon}(v)v \leq 0 \,\,\,\, \mbox{and} \,\,\,\,
\lim_{n \to \infty}\beta_{\epsilon}(v_{n})=\beta(v).
\end{equation}
From (\ref{E30}) and (\ref{E31}), $y=\beta(v) \notin
\frac{1}{\epsilon}\Omega_{+}$, $\Omega_{\epsilon} \subset
A_{\frac{R}{\epsilon},\frac{r}{\epsilon},y}$ and there exists $\tau
\in (0,1]$ such that $\tau v \in
\widehat{\mathcal{M}}_{\epsilon,y}$. Thereby,
$$
\widehat{J}_{\epsilon,y}(\tau v)=E_{\epsilon}(\tau v) \leq
\liminf_{n \to \infty}I_{\epsilon,\lambda_{n},}(\tau v_{n}) \leq
\liminf_{n \to \infty}I_{\epsilon,\lambda_{n}}(v_{n}) \leq
c(\epsilon, B_{r})
$$
which implies
$$
\alpha(R,r,\epsilon,y)\leq c(\epsilon,B_{r}).
$$
On the other hand, since
$$
\alpha(R,r,\epsilon,y)=\alpha(R,r,\epsilon)
$$
we have
$$
\alpha(R,r,\epsilon) \leq c(\epsilon,B_{r}),
$$
contrary to Lemma \ref{L6}, and the proof is complete. $\hfill
\rule{2mm}{2mm}$

\section{Proof of Theorem \ref{T1}}

For $r>0$ and $\epsilon>0$, let $v_{r\epsilon} \in
W^{1,p}_{0}(B_{\frac{r}{\epsilon}}(0))$ be a nonnegative radially
symmetric function such that
$$
E_{\epsilon,B_{r}}(v_{r\epsilon})=c(\epsilon,B_{r})\,\,\, \mbox{and}
\,\,\, E'_{\epsilon,B_{r}}(v_{r\epsilon})=0,
$$
whose existence is proved in \cite[Proposition
4.4]{AlvesFigueiredoSevero1}. For $r>0$ and $\epsilon
>0$, define \linebreak $\Psi_{r}: \frac{1}{\epsilon}\Omega_{-} \to
W_{0}^{1,p}(\Omega_{\epsilon})$ by
$$
\Psi_{r}(y)(x)= \left\{
\begin{array}{l}
v_{r\epsilon}(|x-y|), \,\, x \in B_{\frac{r}{\epsilon}}(y)\\
\mbox{}\\
0, \,\,x \notin B_{\frac{r}{\epsilon}}(y).
\end{array}
\right.
$$
It is immediate that  $\beta_{\epsilon}(\Psi_{r}(y))=y$ for all $y
\in \frac{1}{\epsilon}\Omega_{-}$. In the sequel, we denote by
$I_{\epsilon,\lambda}^{c(\epsilon,B_{r})}$ the set
$$
I_{\epsilon,\lambda}^{c(\epsilon,B_{r})}=\Big\{v \in
\mathcal{M}_{\epsilon,\lambda}:\, I_{\epsilon,\lambda}(v) \leq
c(\epsilon,B_{r}) \Big\}.
$$
We claim that
\begin{equation}\label{E32}
cat I_{\epsilon,\lambda}^{c(\epsilon,B_{r})} \geq cat (\Omega)
\end{equation}
for all $\epsilon \in (0, \epsilon^{*})$ and $\lambda \geq
\lambda^{*}$. In fact, suppose that
$$
I_{\epsilon,\lambda}^{c(\epsilon,B_{r})}=\cup_{i=1}^{n}O_{i}
$$
where $O_{i},i=1,...,n$, is closed and contractible in
$I_{\epsilon,\lambda}^{c(\epsilon,B_{r})}$, that is, there exists
$h_{i} \in C([0,1] \times
O_{i},I_{\epsilon,\lambda}^{c(\epsilon,B_{r})})$ such that, for
every, $v \in O_{i}$,
$$
h_{i}(0,v)=v \,\,\, \mbox{and} \,\,\, h_{i}(1,u)=w_{i}
$$
for some $w_{i}\in I_{\epsilon,\lambda}^{c(\epsilon,B_{r})}$.
Consider
$$
B_{i}=\Psi_{r}^{-1}(O_{i}), \,\,\, i=1,...,n.
$$
The sets $B_{i}$ are closed and
$$
\frac{1}{\epsilon}\Omega_{-}=B_{1}\cup...\cup B_{n}.
$$
Consider the deformation $g_{i}:[0,1] \times B_{i} \to
\frac{1}{\epsilon}\Omega_{+}$ given
$$
g_{i}(t,y)=\beta_{\epsilon}(h_{i}(t,\Psi_{r}(y))).
$$
From Lemma \ref{L7}, the function $g_{i}$ is well defined. Thus,
$B_{i}$ is contractile in $\frac{1}{\epsilon}\Omega_{+}$. Hence,
$$
cat(\Omega)= cat(\Omega_{\epsilon})=
 cat_{\frac{1}{\epsilon} \Omega_{+}}\big(\frac{1}{\epsilon}\Omega_{-}\big)\leq cat I_{\epsilon,\lambda}^{c(\epsilon,B_{r})}
$$
which verifies (\ref{E32}).

Now, we are ready to conclude the proof of Theorem \ref{T1}. From
Proposition \ref{P1} the functional $I_{\epsilon,\lambda}$ satisfies
the Palais-Smale condition provided that $\lambda \geq \lambda^{*}$.
Thus, by Lusternik-Schirelman theory, the functional
$I_{\epsilon,\lambda}$ has at least $cat(\Omega)$ critical points
for all $\epsilon \in (0, \epsilon^{*})$ where $\epsilon^{*}>0$ is
given by Lemma \ref{L6}. The proof is complete. $\hfill
\rule{2mm}{2mm}$


\begin{thebibliography}{99}

%\bibitem{Alves}
%C. O. Alves, P.C. Carri\~ao and E. S. Medeiros, \textit{Multiplicity of solutions for a class of quasilinear problem in exterior
%domains with Newmann conditions}, Abstract and Applied Analisys 03 (2004), 251-268.

%\bibitem{Alves} C. O. Alves, {\it Existence and multiplicity of solution for a class of quasilinear equations},
%Adv. Nonlinear Stud.  5 (2005), no. 1,  73-86.

%\bibitem{Alvesgio1} C. O. Alves and G. M. Figueiredo, \textit{Existence and multiplicity of positive solutions to a p-Laplacian equation in $
%\mathbb{R}^{N}$}, Differential and integral equations  19 (2006), 143-162.

\bibitem{AlvesFigueiredoSevero} C. O. Alves, G. M. Figueiredo and U. B. Severo, \textit{Multiplicity of
positive solutions for a class of quasilinear problems}, Advanced in Differential Equation
14 (2009), 911-942.

\bibitem{AlvesFigueiredoSevero1} C. O. Alves, G. M. Figueiredo and U. B. Severo, \textit{A result of multiplicity of
solutions for a class of quasilinear equations}, to appear in
Proceedings of the Edinburgh Mathematical Society in 2012.

\bibitem{AlvesSoares}C.O. Alves and S.H.M. Soares, \textit{ Multiplicity of positive solutions for a class of nonlinear Schr\"{o}dinger equations}, Advanced in Differential Equations 11 (2010), 1083 - 1102.


\bibitem{AMS1} C.O. Alves, O.H. Miyagaki and S.H.M. Soares, \textit{Multi-bump solutions for a class of quasilinear equations in $\mathbb{R}$}, Communications on Pure and Applied Analysis  11 (2012), 829-844.

\bibitem{AMS2} C.O. Alves, O.H. Miyagaki and S.H.M. Soares, \textit{ On the Existence and Concentration of Positive Solutions to a Class of Quasilinear Elliptic Problems on $\mathbb{R}$}, Mathematische Nachrichten 1 (2011), 1-12.


%\bibitem{Alves10}
%C. O. Alves and M. A. S. Souto, \textit{On existence and concentration behavior of ground state
%solutions for a class of problems with critical growth},
%Comm. Pure and Applied Analysis, 1 (2002), 417-431.
%
\bibitem{MJAlves-Carriao-Miyagaki}
M. J. Alves, P. C. Carri\~ao and  O. H. Miyagaki,  \textit{Soliton solutions to a class of quasilinear
elliptic equations on $\Bbb R$}.  Adv. Nonlinear Stud.  \textbf{7}  (2007), 579--597.
%
%\bibitem{BL} H. Berestycki and P. L. Lions, \textit{Nonlinear scalar
%field equations I: existence of a ground state.} Arch. Rational
%Mech. Anal. \textbf{82}, (1983), 313-346.

\bibitem{BW2} T. Barstch \& Z.Q. Wang,{ \it Existence and multiplicity results for some superlinear
 elliptic problem on $\mathbb{R}^{N}$}, Comm. Partial Differential Equations 20 (1995), no. 9-10, 1725-1741.

\bibitem{BW1} T. Barstch \& Z.Q. Wang,{ \it Multiple positive solutions for a nonlinear Schr\"{o}dinger equation},
  Z. Angew. Math. Phys. 51 (2000), no. 3,  366-384.


\bibitem{Brizhik} L. Brizhik, A. Eremko, B. Piette and W. J. Zakrzewski, \emph{Static
solutions of a $D$-dimensional modified nonlinear Schr{\"o}dinger
equation}, Nonlinearity \textbf{16} (2003) 1481--1497.

\bibitem{Borovskii-Galkin} A. Borovskii and A. Galkin,
\textit{Dynamical modulation of an ultrashort high-intensity laser
pulse in matter.} JETP \textbf{77}, (1983), 562-573.


%\bibitem{brezislieb}
%H. Brezis and E. H. Lieb, A relation between pointwise convergence of functions and convergence
%functionals, Proc. Amer. Math. Soc. 8(1983)486-490.
%
\bibitem{Jeanjean-Colin} M. Colin and L. Jeanjean, \textit{Solutions
for a quasilinear Schr\"{o}dinger equation: a dual approach.} Nonlinear
Anal. \textbf{56}, (2004), 213-226.

\bibitem{Jeanjean-Colin-Squassina} M. Colin and L. Jeanjean and M. Squassina, \textit{Stability and instability results
for standing waves of quasi-linear Schr\"{o}odinger equations.}
Nonlinearity,   \textbf{23}, (2010), 1353-1385.


%
%\bibitem{De Bouard-Hayashi} A. De Bouard, N. Hayashi and J. Saut,
%\textit{Global existence of small solutions to a relativistic
%nonlinear Schrondinger equation.} Comm. Math. Phys. \textbf{189},
%(1997), 73-105.

%\bibitem{costa}
%D. G. Costa. \textit{On a class of elliptic systems in $\mathbb{R}^{N}$,}
% Eletr. J. Diff. Equations. 7 (1994), 1-14.
%
%\bibitem{Djairo}  D. G. de Figueiredo, \textit{Lectures on the Ekeland
%variational principle with applications and detours. Lectures
%notes}, College on Variational Problem in Analysis, Trieste (1988).

%\bibitem{Di}
%E. DiBenedetto
%\textit{$C^{1+\alpha}$ local regularity of weak solutions of
%degenerate results elliptic equations}, Nonl. Analysis TMA 7
%(1983), 827-850.
%
%\bibitem{JMBO e Veve} J. M. do \'O and E. S. de Medeiros,
%\textit{Remarks on least energy solutions for quasilinear elliptic
%problems in $\Bbb R\sp N$.} Electron. J. Differential Equations ,
%\textbf{83}, (2003), 14 pp.


\bibitem{do O-Severo}J.M.B. do \'O and U.B. Severo. \emph{Quasilinear Schr{\"o}dinger equations involving concave and convex
nonlinearities}.  Commun. Pure Appl. Anal.  \textbf{8} (2009),
621--644.


\bibitem{OMS}J.M.B. do \'O, O. H. Miyagaki and S.M.H. Soares,  \textit{Soliton solutions for quasilinear
Schr{\"o}dinger equations: the critical exponential case}, Nonlinear
Anal. \textbf{67}, (2007), 3357-3372.
%
%\bibitem{Ekeland} I. Ekeland, \textit{Convexity methods in
%Hamiltonian Mechanics,} Springer, Berlin, 1990.
%
\bibitem{Floer-Weinstein} A. Floer and A. Weinstein,
\textit{Nonspreading wave pachets for the packets for the cubic
Schrodinger with a bounded potential.} J. Funct. Anal. \textbf{69},
(1986), 397-408.
%
%\bibitem{Gilbarg-Trudinger}D.~Gilbarg and N.~S.~Trudinger, \textit{Elliptic partial
%differential equation of second order}, Second edition,
%Springer-Verlag, Berlin, 1983.
%
%\bibitem{Glowinski-Rappaz} R. Glowinski and J. Rappaz,
%\textit{Approximation of a nonlinear elliptic problem arising in a
%non-Newtonian fluid flow model in glaciology}, Math. Model. Numer.
%Anal. \textbf{37}, (2003),  175-186.

%\bibitem{ghoussoub}
%N. Ghoussoub, \textit{Duality and pertubation methods in critical point
%theory,}Cambridge University Press, Cambridge, 1993


\bibitem{Hartmann} B. Hartmann and W. J. Zakrzewski, \emph{Electrons on hexagonal lattices
and applications to nanotubes}, Phys. Rev. B \textbf{68} (2003)
184302.



\bibitem{Jeanjean-Tanaka} L. Jeanjean and K. Tanaka, \textit{A positive solution for a
nonlinear Schr\"odinger equation on $\mathbb{R}^N$}. Indiana Univ.
Math. \textbf{54} (2005), 443-464.

\bibitem{Kosevich-Ivanov} A. M. Kosevich, B. A. Ivanov and A. S.
Kovalev, \textit{Magnetic solitons in superfluid films.} J. Phys.
Soc. Japan \textbf{50}, (1981), 3262-3267.

%\bibitem{livro do Kavian} O. Kavian, \textit{Introduction \'a la th\'eorie des points critiques
%et applications aux probl\`emes elliptiques}, Springer-Verlag,
%Paris, 1993.

%\bibitem{Gongbao}
%Li Gongbao, \textit{Some properties of weak solutions of nonlinear
%scalar field equations},
%Annales Acad. Sci. Fenincae, series A. 14 (1989), 27-36.
%

\bibitem{Kurihura} S. Kurihura, \textit{Large-amplitude
quasi-solitons in superfluids films.} J. Phys. Soc. Japan
\textbf{50}, (1981), 3262-3267.
%
%\bibitem{Laedke-Spatschek} E. Laedke and K. Spatschek,
%\textit{Evolution theorem  for a class of perturbed envelope soliton
%solutions.} J. Math. Phys. \textbf{24}, (1963), 2764-2769.
%
%\bibitem{Lions}  P. L. Lions, \textit{The concentration-compacteness principle in the
%calculus of variations. The locally compact case, Part II}, Ann.
%Inst. H. Poincar\'{e} Anal. Non Lin\'{e}are, \textbf{1}, (1984),
%223-283.
%
\bibitem{Liu-Wang I} J. Liu and Z. Q. Wang, \textit{Soliton
solutions for quasilinear Schr\"{o}dinger equations I}, Proc. Amer.
Math. Soc. \textbf{131}, 2, (2002), 441--448.
%
\bibitem{Liu-Wang II} J. Liu, Y. Wang and Z. Wang, \textit{Soliton
solutions for quasilinear Schr�dinger equations II.} J. Differential
Equations \textbf{187}, (2003), 473-493.
%
\bibitem{Liu-wang-wang} J. Liu, Y. Wang and Z. Q. Wang, \textit{Solutions
for Quasilinear Schr\"odinger Equations via the Nehari Method},
Comm. Partial Differential Equations, \textbf{29}, (2004) 879--901.
%
\bibitem{Makhankov-Fedyanin} V. G. Makhankov and V. K. Fedyanin,
\textit{Non-linear effects in quasi-one-dimensional models of
condensed matter theory.} Phys. Reports \textbf{104}, (1984), 1-86.

%\bibitem{Moser}
%J. Moser, \textit{A new proof of de Giorgi's theorem concerning the
%regularity problem for elliptic differential equations,}
%Comm. Pure Apll. Math. 13 (1960), 457-468.

%
\bibitem{Poppenberg-Schmitt-Wang} M. Poppenberg, K. Schmitt and Z.
Q. Wang, \textit{On the existence of soliton solutions to
quasilinear Schr\"odinger equations}, Calc. Var. Partial
Differential Equations \textbf{14}, (2002), 329--344.
%

\bibitem{Rabinowitz} P. H. Rabinowitz,{ \it On a class of nonlinear Schr\"{o}dinger equations},
Z. Angew Math. Phys. 43 (1992), no. 2, 270-291.

\bibitem{Ritchie} B. Ritchie, \textit{Relativistic self-focusing and
channel formation in laser-plasma interactions.} Phys. Rev. E
\textbf{50}, (1994), 687-689.
%
%\bibitem{Serrin} J. Serrin, \textit{Local behavior of solutions of quasi-linear
%equations}, Acta Math. 111 (1964), 247-302.
%
\bibitem{severo}
U.B. Severo, \textit{Estudo de uma classe de equa\c{c}\~{o}es de Schr\"{o}dinger quase-lineares}.
Doct. dissertation, Unicamp, 2007.
%
\bibitem{severo2}
 U.B. Severo, \textit{Existence of weak solutions for quasilinear elliptic equations involving the p-Laplacian}. Electron. J. Differential Equations (2008), no. 56, 1-16.


\bibitem{Elves1} E.A.B. Silva  and G.F. Vieira, \textit{Quasilinear asymptotically periodic Schr\"odinger
equations with critical growth}, Calc. Var. 39 (2010), 1�33

\bibitem{Elves2} E.A.B. Silva and G.F. Vieira, \textit{ Quasilinear asymptotically periodic Schr\"odinger equations with subcritical
growth,} Nonlinear Anal 72 (2010), 2935 - 2945.

\bibitem{Takeno-Homma} S. Takeno and S. Homma, \textit{Classical
planar Heinsenberg ferromagnet, complex scalar fields and nonlinear
excitations.} Progr. Theoret. Physics \textbf{65}, (1981), 172-189.
%
%\bibitem{tolksdorf} P. Tolksdorff, \textit{Regularity for a more general class of quasilinear
%elliptic equations}, J. Differential Equations, \textbf{51} (1984),
%126-150.
%
%\bibitem{Trudinger} N. S. Trudinger, \textit{On Harnack type inequalities and their
%application to quasilinear elliptic equations}, Comm. Pure Appl.
%Math. 20 (1967) 721-747.

%\bibitem{Willem}
%M. Willem , \textit{\em Minimax Theorems}, Birkh�user, 1996.

\end{thebibliography}
\end{document}